\documentclass[12pt]{article}

\title{Sequential robust efficient estimation for nonparametric 
autoregressive models}

\author{Ouerdia Arkoun
\thanks{
Normandie Universit\'e, Universit\'e de Rouen, 
Laboratoire de Math\'ematiques Rapha\"el Salem,
CNRS, UMR 6085, Avenue de l'universit\'e, BP 12, 
76801 Saint-Etienne du Rouvray Cedex, 
France.
\newline
email:  Ouerdia.Arkoun@gmail.com}
\and Serguei Pergamenchtchikov
\thanks{
Normandie Universit\'e, Universit\'e de Rouen, 
Laboratoire de Math\'ematiques Rapha\"el Salem,
CNRS, UMR 6085, Avenue de l'universit\'e, BP 12, 
76801 Saint-Etienne du Rouvray Cedex, 
France and
Department of Mathematics and Mechanics, Tomsk State University, Lenin str. 36, 
634041 Tomsk, Russia
\newline
email:  Serge.Pergamenchtchikov@univ-rouen.fr}
\date{}
}
\usepackage{amssymb}
\usepackage{amsfonts}
\usepackage{amsmath}
\usepackage{amsthm}
\newtheorem{theorem}{Theorem}[section]
\newtheorem{proposition}[theorem]{Proposition}
\newtheorem{lemma}[theorem]{Lemma}

\newtheorem{remark}[theorem]{Remark}

\newcommand\e{\varepsilon}
\def\bbr{{\mathbb R}}

\newcommand{\wh}{\widehat}
\newcommand{\wt}{\widetilde}

\newcommand\cH{{\cal H}}
\newcommand\cG{{\cal G}}
\newcommand\cF{{\cal F}}

\newcommand\cN{{\cal N}}
\newcommand\cP{{\cal P}}

\newcommand\cR{{\cal R}}

\newcommand\cU{{\cal U}}

\def\text#1{\hbox{#1}}
\def\proof{{\noindent \bf Proof. }}

\def\endproof{\mbox{\hfill $\qed$}}

\def\E{{\bf E}}
\def\P{{\bf P}}

\def\M{{\bf M}}
\def\r{{\bf r}}
\def\m{{\bf m}}
\def\C{{\bf C}}
\def\L{{\bf L}}

\def\b{{\bf b}}

\def\s{{\bf s}}
\def\R{{\bf R}}
\def\Chi{{\bf 1}}
\def\d{\mbox{d}}
\def\build #1_#2{\mathrel{\mathop{\kern 0pt #1}\limits_\zs{#2}}} 
\newcommand{\zs}[1]{{\mathchoice{#1}{#1}{\lower.25ex\hbox{$\scriptstyle#1$}}
{\lower0.25ex\hbox{$\scriptscriptstyle#1$}}}}

\numberwithin{equation}{section}

\begin{document}

\maketitle

\begin{abstract}

We construct efficient robust truncated sequential  estimators
for the pointwise estimation problem  
in nonparametric autoregression models with smooth coefficients.
For Gaussian models we propose
 an adaptive procedure based on the constructed sequential estimators. 
 The minimax nonadaptive and adaptive convergence rates are established. It turns out that 
in this case these rates are the same  as for regression models.
\end{abstract}

\par
{\bf AMS (2000) Subject Classification :}
primary 62G07,62G08; secondary 62G20.

\section{Introduction}\label{sec:In}

One of the standard linear models in the general theory of time series is the 
autoregressive model (see, for example, \cite{And70} and the references therein).
Natural extensions for such models are nonparametric autoregressive models which are defined
by
\begin{equation}\label{sec:In.1}
y_\zs{k} = S(x_\zs{k})y_\zs{k-1}+ \xi_\zs{k}\,,\quad 1\le k\le n\,,
\end{equation}
where $S(\cdot)$ is unknown function, the design
 $x_\zs{k}=k/n$ and 
 the noise $(\xi_\zs{k})_\zs{1\le k\le n}$ are i.i.d. 
unobservable centered random variables, i.e. $\E\,\xi_\zs{1}=0$.

It should be noted that the varying coefficient
principle is well known in the regression analysis.
It permits the use of a more complex forms for  regression coefficients
and, therefore, the models constructed via this method
are more adequate for applications (see, for instance, 
\cite{FaZh}, \cite{LuYaZh}). In this paper we consider 
the varying coefficient autoregressive models \eqref{sec:In.1}. 
There is a number of papers which consider these  models 
such as  \cite{Da96}, \cite{Da96-1} and \cite{Be00}.
In all these papers, the authors propose some asymptotic (as $n\to\infty$) methods  for different 
identification studies without considering optimal estimation issues.
To our knowledge, for the first time the minimax estimation 
problem
for the model \eqref{sec:In.1} 
has been treated in  
\cite{ArPe08} in the nonadaptive case, 
i.e. for the known regularity of the function
$S$. 
Then, in  \cite{Ark11} it is proposed to use the 
sequential analysis method  
for the adaptive  pointwise estimation  problem 
in the case when the unknown H\"older regularity is less than one, i.e when the function $S$
is not differentiable.
It turns out that it is only the sequential analysis
  method allows to construct
 an adaptive pointwise estimation procedure
for the models \eqref{sec:In.1}.
That is why, in this paper, we study sequential estimation methods  
 for the smooth 
function $S$. We consider the pointwise estimation at a fixed point 
$z_\zs{0}\in]0;1[$ in the two cases: when 
the H\"older regularity is known  and when
the H\"older regularity is unknown, i.e.
the adaptive estimation. In the first case we consider this problem in the robust 
setting, i.e. we assume that the distribution of the random variables 
$(\xi_\zs{j})_\zs{j\ge 1}$
in \eqref{sec:In.1}
belongs to some functional class and we consider the estimation problem with respect 
to robust risks which have an additional supremum over all distributions from
some fixed class. For nonparametric regression models such risks was introduced 
in \cite{GaPe06}
for the pointwise estimation problem and in \cite{GaPe09} for the quadratic risks.
Later, for the quadratic risks the same approach was used in \cite{KoPe12} 
for regression model in continuous time. Motivated by this facts, we consider the adaptive estimation problem for the Gaussian models \eqref{sec:In.1}.
More precisely, we assume that the function $S$ belongs to a H\"{o}lder class with some
unknown regularity $1 < \beta \leq 2$. Unfortunately, we can not use directly
the sequential procedure from
\cite{Ark11} for the adaptive estimation of such functions.  Since to obtain an optimal 
rate for the function with  $\beta>1$ we have to
take into account  the Taylor expansion of the function $S$ at $z_\zs{0}$ of the order $1$.
To study the Taylor expansion for sequential procedures one needs to control
the  behavior of the stopping time. Indeed, one needs to keep the stopping time near of the number of observations.
This can not be done by the procedure
from \cite{Ark11} since one needs to use the unknown function
 $S$.
In this paper we construct a sequential adaptive estimate for smooth functions
and we find an adaptive minimax convergence rate  for smooth functions.
In Section~\ref{sec:Sp}, we present the standard notations used in sequel of the paper. 
We describe in detail the statement of the problem and main results in Section~\ref{sec:Mr}. 
In Section~\ref{sec:Pr}, we will study some properties of kernel estimators
for the model \eqref{sec:In.1}  and in Section~\ref{sec: Stp}, 
we study  properties of  stopping time for constructed sequential 
procedure. Section~\ref{sec:Prf} is devoted to the asymptotic 
upper bound and lower bound for the risk of the sequential kernel estimators.
In Section~\ref{sec:Nur} we illustrate the obtained results by numerical examples. 
Finally, we give the appendix which contains some technical results. 

\vspace{2mm}

\section{Sequential procedures}\label{sec:Sp}
\subsection{Main conditions}
We assume that  in the model \eqref{sec:In.1} the i.i.d. random
variables  $(\xi_\zs{k})_\zs{1\le k\le n}$  have a density  $p$
(with respect to the Lebesgue measure)  from the functional class $\cP_\zs{\varsigma}$
defined as
\begin{align}\nonumber
\cP_\zs{\varsigma}:=&
\left\{
p\ge 0\,:\,\int^{+\infty}_\zs{-\infty}\,p(x)\,\d x=1\,,\quad
\int^{+\infty}_\zs{-\infty}\,x\,p(x)\,\d x= 0 \,,\right.
\\[2mm] \label{2.1} 
&\quad\left.
\int^{+\infty}_\zs{-\infty}\,x^2\,p(x)\,\d x= 1
\quad\mbox{and}\quad
\sup_\zs{k\ge 1}\,\frac{1}{\varsigma^{k}(2k-1)!!}
 \int^{+\infty}_\zs{-\infty}\,|x|^{2k} \,p(x)\,\d x\le 1
\,
\right\}\,,
\end{align}
where $\varsigma\ge 1$ is some fixed parameter.
  Note that the $(0,1)$-Gaussian density 
belongs to $\cP$. In the sequel we denote this density by $ p_\zs{0}$.
It is clear that  
for any $q>0$
\begin{equation}\label{sec:Sp.0} 
\s^{*}_\zs{q}=\sup_\zs{p\in\cP_\zs{\varsigma}}\,\E_\zs{p}\,|\xi_\zs{1}|^{q}
<\infty\,,
\end{equation}
where $\E_\zs{p}$ is the expectation with respect to the density $p$ from
$\cP_\zs{\varsigma}$.
\noindent  To obtain the stable (uniformly with respect to the function $S$ ) model \eqref{sec:In.1}, we assume that for some fixed $0<\e<1$ and $L>0$ the unknown function $S$ belongs to 
the $\varepsilon$ - {\em stability set} 
\begin{equation}\label{sec:Sp.1}
\Theta_\zs{\e,L} = \left\{S\in \C_\zs{1}([0,1],\bbr) : \|S\| \le 1-\e
\quad\mbox{and}\quad
\|\dot S\|\le L 
 \right\}\,, 
\end{equation}
where 
$\C_\zs{1}[0,1]$ is the Banach space of continuously differentiable
$[0,1]\to\bbr$ functions and $\|S\| = \sup_\zs{0\le x \leq1}|S(x)|$. 
 Similarly to \cite{GaPe06} and \cite{ArPe08} we make use  of the family 
of the {\em weak stable local H\"older classes} at the point $z_\zs{0}$
\begin {equation}\label{sec:Sp.1-0}
  \cU^{(\beta)}_\zs{n}(\e,L,\epsilon^{*}_\zs{n}) =
\left\{S\in \Theta_\zs{\e,L}\,:\,
|\Omega_\zs{h}(z_\zs{0},S)| \le \epsilon^{*}_\zs{n} h^{\beta} 
\right\}\,,
\end{equation}
where 
$$
\Omega_\zs{h}(z_\zs{0},S) =  \int_\zs{-1}^{1}(S(z_\zs{0}+uh)- S(z_\zs{0}))\,\d u
$$ 
and $\beta=1+\alpha$ is the regularity 
parameter with
$0<\alpha<1$. Moreover, we assume that the {\em weak H\"older constant}
$\epsilon^{*}_\zs{n}$ goes to zero, i.e. $\epsilon^{*}_\zs{n}\to 0$ as $n\to\infty$. 
Moreover, we define
the corresponding {\em strong stable  local  H\"older class} at the point $z_\zs{0}$
as

\begin{equation}\label{sec:Sp.2}
 \cH^{(\beta)}(\e,L,L^{*}) 
= \left\{S \in \Theta_\zs{\e,L}\,:\,\Omega^{*}(z_\zs{0},S) \le L^{*} \right\}\,,
\end{equation}
where
$$
\Omega^*(z_\zs{0},S) = \sup_\zs{x\in [0,1]} 
\frac{|\dot{S}(x) - \dot{S}(z_\zs{0})|}{|x-z_\zs{0}|^\alpha}\,.
$$
We assume that the regularity $\underline{\beta}\le \beta \le \overline{\beta}$, where
$\underline{\beta}=1+\underline{\alpha}$ and $\overline{\beta}=1+\overline{\alpha}$ for some fixed parameters
$0<\underline{\alpha}<\overline{\alpha}<1$.

\begin{remark}\label{Re.sec:Sp.1}
Note that for the regression models
the weak H\"older class was introduced in \cite{GaPe06}
for the efficient pointwise estimation. It is clear that it is more large than
usual one with the same H\"older constant, i.e.
$$
\cH^{(\beta)}(\e,L,\epsilon^{*}_\zs{n}) \subseteq 
\cU^{(\beta)}_\zs{n}(\e,L,\epsilon^{*}_\zs{n})\,.
$$  
It should be noted also that for  diffusion processes the local weak H\"older 
class was used in \cite{GaPe05} and \cite{GaPe06-1} for sequential and truncated
sequential efficient pointwise estimation respectively. Moreover,  in
\cite{GaPe11} these sequential pintwise efficient estimators were used to construct 
 adaptive efficient model selection procedures in $\L_\zs{2}$ for diffusion processes. 
\end{remark}

\subsection{Nonadaptive procedure}

First,  we study a nonadaptive estimation problem 
for the function $S$ from the functional class \eqref{sec:Sp.1-0}
of the known regularity $\beta=1+\alpha$. As we will see later to construct an efficient sequential procedure we need to use  $S$ 
as a procedure parameter. So we propose to use the first $\nu$ observations for the auxiliary
estimation of $S(z_\zs{0})$. In this step we use usual kernel estimate, i.e. 
\begin{equation}\label{sec:Sp.3}
\wh{S}_\zs{\nu}=\frac{1}{A_\zs{\zs{\nu}}}\,
\sum^{\nu}_\zs{j=1}\,Q(u_\zs{j})\,y_\zs{j-1}\,y_\zs{j}\,,\quad
A_\zs{\nu}=\sum^{\nu}_\zs{j=1}\,Q(u_\zs{j})\,y_\zs{j-1}^2\,,
\end{equation}
where  the kernel $Q(\cdot)$ is the indicator function of the interval $[-1;1]$;
$u_\zs{j} = (x_\zs{j}-z_\zs{0})/h$ and $h$ is some positive 
bandwidth. 
In the sequel for any $0\le k< m\le n$ we set
\begin{equation}\label{sec:Sp.3-0}
A_\zs{k,m}=\sum^{m}_\zs{j=k+1}\,Q(u_\zs{j})\,y_\zs{j-1}^2\,,
\end{equation}
i.e. $A_\zs{\nu}=A_\zs{0,\nu}$.
It is clear that to estimate $S(z_\zs{0})$
on the basis of the kernel estimate with the
kernel $Q$ we can use the observations 
$(y_\zs{j})_\zs{k_\zs{*}\le j\le k^{*}}$, where
\begin{equation}\label{sec:Sp.7}
k_\zs{*}= [nz_\zs{0} - nh ] + 1
\quad\mbox{and}\quad
k^{*}=  [nz_\zs{0} + nh]\,.
\end{equation}
Here $[a]$ is the integral part of a number $a$.
So for the first estimation we chose $\nu$ as 
\begin{equation}\label{sec:Sp.9-0}
\nu=\nu(h,\alpha)=k_\zs{*}+\iota\,,
\end{equation}
where 
$$
\iota = \iota(h,\alpha) = [\wt{\epsilon}nh]+1 
\quad\mbox{and}\quad
\wt{\epsilon}=
\wt{\epsilon}(h,\alpha)=h^{\alpha}/\ln n
\,.
$$

Next, similarly to \cite{Ark11}, 
we use a some kernel  sequential procedure 
 based on the observations
$(y_\zs{j})_\zs{\nu\le j\le n}$.
To transform the kernel estimator  in the linear function of observations
we replace the number of observations $n$ by the following stopping time
\begin{equation}\label{sec:Sp.4}
 \tau_\zs{H} = \inf\{ k\geq \nu+1: A_\zs{\nu,k} \ge H \},
 \end{equation}
where $\{\emptyset\}=n$ and the positive threshold $H$ will be chosen as a positive 
random variable measurable with respect to the $\sigma$ - field 
$\{y_\zs{1},\ldots,y_\zs{\nu}\}$. Therefore, we get 

\begin{equation}\label{sec:Sp.5}
S^{*}_\zs{h}=\frac{1}{H}\,
\left(\sum^{\tau_\zs{H}-1}_\zs{j=\nu+1}\,Q(u_\zs{j})\,y_\zs{j-1}\,y_\zs{j}\,+\,\varkappa_\zs{H}\,Q(u_\zs{\tau_\zs{H}})\,y_\zs{\tau_\zs{H}-1}\,y_\zs{\tau_\zs{H}}\right)  
\Chi_\zs{(A_\zs{\nu,n}\geq H)}\,,
\end{equation} 
where the correcting  coefficient $\varkappa_\zs{H}$ on the set 
$\{A_\zs{\nu,\tau_\zs{H}}\geq H\}$ is defined as
$$
 A_\zs{\nu,\tau_\zs{H}-1} + \varkappa_\zs{H}\,Q(u_\zs{\tau_\zs{H}})\,
y^2_\zs{\tau_\zs{H}-1} = H
$$
and $\varkappa_\zs{H} =1$  on the set $\{A_\zs{\nu,\tau_\zs{H}}< H\}$.

Now, to obtain an efficient estimate we need to use the all $n$ observations, i.e. asymptotically 
for sufficiently large $n$ the stopping time $\tau_\zs{H} \approx n$. 
Similarly to  \cite{KoPe84}, one can show that $\tau_\zs{H}\approx \gamma(S)\,H$ as $H\to\infty$, where
\begin{equation}\label{sec:Sp.6}
\gamma(S)=1-S^2(z_\zs{0})\,.
\end{equation} 
Therefore, to use asymptotically all  observations we have to chose 
$H$ as the number observations divided by $\gamma(S)$. But in our case 
we use $k^{*} - k_\zs{*}$ observations to estimate $S(z_\zs{0})$, 
Therefore, to obtain optimal estimate we need to define $H$ as $(k^{*} - k_\zs{*})/\gamma(S)$.  
taking into account that $ k^{*} - k_\zs{*}\approx 2nh$ and that $\gamma(S)$ is unknown
we define the threshold $H$ as

\begin{equation}\label{sec:Sp.8}
H=H(h,\alpha)=\phi nh\,,\quad 
\phi=\phi(h,\alpha)=\frac{2(1-\wt{\epsilon})}{\gamma(\wt{S}_\zs{\nu})}
\,,
\end{equation}
where $\wt{S}_\zs{\nu}$ is the projection of the estimator 
$\wh{S}_\zs{\nu}$ in the interval $]-1+\e,1-\e[$, i.e.
$$
\wt{S}_\zs{\nu} = \min (\max (\wh{S}_\zs{\nu},-1+\e),1-\e)\,.
$$
In this paper we chose the bandwidth $h$ in the following form
\begin{equation}\label{sec:Sp.9}
h=h(\beta)=(\kappa_\zs{n})^{\frac{1}{2\beta+1}}\,,
\end{equation}
where the sequence $\kappa_\zs{n}>0$ such that
\begin{equation}\label{sec:Sp.10}
\kappa_\zs{*}=\liminf_\zs{n\to\infty} n\kappa_\zs{n}>0
\quad\mbox{and}\quad
\lim_\zs{n\to\infty} n^{\delta}\,\kappa_\zs{n}=0
\end{equation}
for any $0<\delta<1$.

\medskip

\subsection{Adaptive procedure}
We will construct 
an adaptive minimax sequential estimation  
for the function $S$ from the functional class \eqref{sec:Sp.2}
of the unknown regularity $\beta$. To this end we will use the modification of the adaptive Lepskii
method proposed in \cite{Ark11} based on the sequential estimators \eqref{sec:Sp.5}.
We set
\begin{equation}\label{sec:Ad.1}
d_\zs{n} = \frac{n}{\ln n}
 \quad\mbox{and}\quad 
N(\beta)=\left(d_\zs{n}\right)^{\frac{\beta}{2\beta+1}}
\,. 
\end{equation}
Moreover, we chose the bandwidth $h$ in the form \eqref{sec:Sp.9}
with $\kappa_\zs{n}=1/d_\zs{n}$, i.e. we set
\begin{equation}\label{sec:Ad.1-1}
\check{h}=
\check{h}(\beta)=\left(\frac{1}{d_\zs{n}}\right)^{\frac{1}{2\beta+1}}\,.
\end{equation}
We define the grids on the intervals 
$[\underline{\beta}\,,\,\overline{\beta}]$
and $[\underline{\alpha}\,,\,\overline{\alpha}]$
 as
\begin{equation}\label{sec:Ad.2}
\beta_\zs{k}=\underline{\beta}+\frac{k}{m}(\overline{\beta}-\underline{\beta})
\quad\mbox{and}\quad
\alpha_\zs{k}=\underline{\alpha}+\frac{k}{m}
(\overline{\alpha}-\underline{\alpha})
\end{equation}
for $0\le k\le m$ with $m=[\ln d_\zs{n}]+1$, and we set 
$$
N_\zs{k}=N(\beta_\zs{k})
\quad\mbox{and}\quad 
 \check{h}_\zs{k}=\check{h}(\beta_\zs{k})\,.
$$
Replacing in \eqref{sec:Sp.9-0} and 
\eqref{sec:Sp.8} the parameters $h$ and $\alpha$ we define
$$
\check{\nu}_\zs{k}=\nu(\check{h}_\zs{k},\alpha_\zs{k})
\quad\mbox{and}\quad
\check{H}_\zs{k}=H(\check{h}_\zs{k},\alpha_\zs{k})
\,.
$$
Now 
using these parameters 
in the estimators
\eqref{sec:Sp.3} and \eqref{sec:Sp.5}
we set $
\check{S}_\zs{k}=S^{*}_\zs{\check{h}_\zs{k}}(z_\zs{0})$
and
\begin{equation}\label{sec:Ad.3}
\check{\omega}_\zs{k}=\max_\zs{0\le j\le k}\,
\left(
|\check{S}_\zs{j}
-
\check{S}_\zs{k}
|
-
\frac{\check{\lambda}}{N_\zs{j}}
\right)\,,
\end{equation}
where 
$$
\check{\lambda}>\check{\lambda}_\zs{*}=4\sqrt{2}
\left(
\frac{\overline{\beta}-\underline{\beta}}{(2\underline{\beta}+1)(2\overline{\beta}+1)}
\right)^{1/2}\,.
$$
In particular, if $\underline{\beta}=1$ and $\overline{\beta}=2$
we get $\check{\lambda}_\zs{*}=4(2/15)^{1/2}$.
\noindent 
We also define the 
optimal index as
\begin{equation}\label{sec:Ad.4}
\check{k}=\max\left\{0\le k\le m\,:\,\check{\omega}_\zs{k}\le
\frac{\check{\lambda}}{N_\zs{k}}\right\}\,.
\end{equation}

\noindent The adaptive estimator is now defined as
\begin{equation}\label{sec:Ad.5}
\wh{S}_\zs{a,n}=S^{*}_\zs{\check{h}_k}
\quad\mbox{and}\quad
\check{h}_k = \check{h}_\zs{\check{k}}\,.
\end{equation}

\begin{remark}\label{Re.sec:Sp.2}
It should be noted that in the difference from the
usual adaptive pointwise estimation (see, for example, \cite{Le90}, \cite{GaPe01},
\cite{Ark11} and al.)  the threshold 
$\check{\lambda}$ in  \eqref{sec:Ad.3}
does not depend 
on the  parameters $L>0$ and $L^{*}>0$ 
of the H\"older class \eqref{sec:Sp.2}.
\end{remark}

\medskip

\section{Main results}\label{sec:Mr}

\subsection{Robust efficient estimation}

The problem is to  estimate the function $S(\cdot)$ 
at a fixed point $z_\zs{0}\in ] 0, 1[$, i.e. the value $S(z_\zs{0})$.
For this problem we make use of the risk proposed in \cite{ArPe08}. Namely,
  for any estimate $\wt{S}=\wt{S}_\zs{n}(z_\zs{0})$ 
(i.e. any measurable with respect to the observations $(y_\zs{k})_\zs{1\le k\le n}$ function)
we define the following robust risk
\begin{equation}\label{sec:Mr.1}
\cR_\zs{n}(\wt{S}_n,S)= \sup_\zs{p\in \cP_\zs{\varsigma}}\,
\E_\zs{S,p}|\wt{S}_\zs{n}(z_\zs{0})-S(z_\zs{0})|\,,
\end{equation}
where $\E_\zs{S,p}$ is the expectation taken with respect to the 
distribution $\P_\zs{S,p}$ of the vector $(y_1,...,y_\zs{n_0})$ in 
\eqref{sec:In.1} corresponding to the function S and the density 
$p$ from $\cP_\zs{\varsigma}$. 

With the help of the function $\gamma(S)$ defined in \eqref{sec:Sp.6}, we 
describe the sharp lower bound for the minimax risks 
with the normalizing coefficient 
\begin{equation}\label{sec:Mr.2}
\varphi_\zs{n}=n^{\frac{\beta}{2\beta+1}}\,.
\end{equation}

\begin{theorem}\label{Th.sec:Mr.1}
For any $0<\e<1$
 \begin{equation}\label{sec:Mr.3}
\underline{\lim}_\zs{n\to\infty} \, \inf_\zs{\wt{S}}\,
\sup_\zs{S\in  \cU^{(\beta)}_\zs{n}(\e,L,\epsilon^{*}_\zs{n})} \gamma^{-1/2}(S)\,
 \varphi_n\,
\cR_\zs{n}(\wt{S}_\zs{n},S) \ge
 \E |\eta|\,,
\end{equation}
where $\eta$ is a Gaussian random variable with the parameters $(0,1/2)$.
\end{theorem}

\vspace{10mm}
Now we give the upper bound for the minimax risk of the sequential kernel estimator defined in \eqref{sec:Sp.5}.

\vspace{2mm}
\begin{theorem}\label{Th.sec:Mr.2}
The estimator \eqref{sec:Sp.3} with the parameters 
\eqref{sec:Sp.8} -- \eqref{sec:Sp.9} and $\kappa_\zs{n}=n^{-1}$
 satisfies the following inequality 
$$
\overline{\lim}_\zs{n\to\infty}\,
\sup_\zs{S\in  \cU^{(\beta)}_\zs{n}(\e,L,\epsilon^{*}_\zs{n})} \gamma^{-1/2}(S)\,
\varphi_n \,
\cR_\zs{n}(\wh{S}_\zs{a,n},S) \le \E |\eta|\,,
 $$
where $\eta$ is a Gaussian random variable with the parameters $(0,1/2)$. 
\end{theorem}

\begin{remark}\label{Re.sec:Mr.1}
Theorems \ref{Th.sec:Mr.1} and \ref{Th.sec:Mr.1} imply that the  estimator \eqref{sec:Sp.3}, with the parameters \eqref{sec:Sp.9}
 is asymptotically robust efficient with respect to class $\cP_\zs{\varsigma}$.
\end{remark}

\subsection{Adaptive estimation}

Now we consider the Gaussian model \eqref{sec:In.1}, i.e. assume that the random variables
$(\xi_\zs{j})_\zs{j\ge 1}$ are  $\cN(0,1)$. 
 The problem is to estimate the function $S$ at a fixed point 
$z_\zs{0}\in ]0,1[,$ i.e. the value $S(z_\zs{0})$. For any estimate 
$\tilde{S}_\zs{n}$ of $S(z_\zs{0})$ 
(i.e.  any measurable with respect to the observations 
$(y_\zs{k})_\zs{1\le k\le n}$ function), we define
the adaptive risk for the functions $S$ from $\cH^{(\beta)}(\e,L,L^{*})$ as 
\begin{equation}\label{sec:Mr.4}
\cR_\zs{a,n}(\tilde{S}_n)=
\sup_\zs{\beta\in [\underline{\beta};\overline{\beta}]}\,
\sup_\zs{S\in \cH^{(\beta)}(\e,L,L^{*})}\,N(\beta)\,
\E_\zs{S}|\tilde{S}_\zs{n}-S(z_\zs{0})|\,,
\end{equation}
where $ N(\beta)$ is defined in
\eqref{sec:Ad.1}, 
$\E_\zs{S}=\E_\zs{S,p_\zs{0}}$ is the expectation taken with respect to the distribution 
$\P_\zs{S}=\P_\zs{S,p_\zs{0}}$.

First we give the lower bound for the minimax risk. 
We show that with the convergence rate $N(\beta)$ the lower bound for the minimax risk is strictly positive.

Now we give the upper bound for the minimax risk of the sequential kernel estimator defined in \eqref{sec:Sp.5}.

First we give the lower bound for the minimax risk. We show that with the convergence rate $N(\beta)$ the lower bound for the minimax risk is strictly positive.
\begin{theorem}\label{Th.6.1}
There exists $L^{*}_\zs{0} > 0$ such that for all  $L^{*}> L^{*}_\zs{0}$, the risk \eqref{sec:Mr.4} 
admits the following lower bound: 
\begin{equation*}
\liminf_\zs{n \to \infty}\,\inf_\zs{\tilde{S}_n}\,\cR_\zs{a,n}(\tilde{S}_n)>0\,,
\end{equation*}
where the infimum is taken over all estimators $\tilde{S_n}.$
\end{theorem}
\noindent 
The proof of this theorem is given in \cite{Ark11}.

To obtain an upper bound for the adaptive risk \eqref{sec:Mr.4}
of  the procedure \eqref{sec:Ad.5}
 we need to study 
 the family  
$(S^{*}_\zs{h})_\zs{\underline{\alpha} \leq\alpha\leq \overline{\alpha}}$.

\begin{theorem}\label{Th.sec:Mr.1_1} The sequential procedure  \eqref{sec:Sp.5}
 with the bandwidth $h$ defined in \eqref{sec:Sp.9}
for $\kappa_\zs{n}=\ln n/n$
 satisfies the following 
property
$$
\limsup_\zs{n\to \infty}
\sup_\zs{\underline{\alpha} \leq\alpha\leq \overline{\alpha}}\,
(\Upsilon_n(h))^{-1}\,
\sup_\zs{S\in\cH^{(\beta)}(\e,L,L^{*})}\,\sup_\zs{p\in \cP_\zs{\varsigma}}\,
\E_\zs{S,p}|S^{*}_\zs{h} - S(z_\zs{0}))| < \infty $$
where $\Upsilon_n(h) = h^{\beta} + (nh)^{-1/2}$.
\end{theorem}

\noindent Using this theorem we can establish the minimax property
for the procedure \eqref{sec:Mr.4}.

\begin{theorem}\label{Th.6.2}
 The estimation procedure \eqref{sec:Mr.4} satisfies the following
asymptotic property
\begin{equation}\label{4.6-0}
\limsup_\zs{n\to\infty}\,\cR_\zs{a,n}(\wh{S}_\zs{a,n}) < \infty\,. 
\end{equation}
\end{theorem}

\begin{remark}\label{Re.sec:Mr.2}
Theorem~\ref{Th.6.1} gives the lower bound for the adaptive risk, i.e. the convergence rate $N(\beta)$ is best for the adapted risk. Moreover, by Theorem~\ref{Th.6.2} the adaptive estimates \eqref{sec:Ad.5} possesses this convergence rate. In this case, this estimates is called optimal in sense of the adaptive risk \eqref{sec:Mr.4}
\end{remark}


\section{Properties of $\wh{S}_\zs{\nu}$}\label{sec:Pr}

We start with studying the properties of the estimate \eqref{sec:Sp.8}.
 To this end for any  $q> 1$  we set
\begin{equation}\label{sec:Pr.1}
\varrho^{*}_\zs{q}=
 \frac{(12(1+\kappa_\zs{*}))^{q}}{(\varepsilon\kappa_\zs{*})^{2q}}
\left(\b^{*}_\zs{q}\left(r^{*}_\zs{q}\s^{*}_\zs{q}+\s^{*}_\zs{2q}+1\right) +
2(1+L)^{q}\, r^{*}_\zs{2q}
\right)
\,,
\end{equation}
where 
$$
\r^{*}_\zs{q}= 2^{q-1}\left( |y_\zs{0}|^{q}+
\s^{*}_\zs{q}\,\left(\frac{1}{\e}\right)^{q}
\right)
\quad\mbox{and}\quad
\b^{*}_\zs{q} = \frac{18^{q} q^{3q/2}}{(q-1)^{q/2}}\,.
$$
Now we obtain a non asymptotic upper bound for the tail probability
for the deviation
\begin{equation}\label{sec:Pr.1-1}
\wh{\Delta}_\zs{\nu}=
\wh{S}_\zs{\nu}-S(z_\zs{0})\,.
\end{equation}

\medskip

\begin{proposition}\label{Le.S.1}
For any  $q>1$, $h>0$ and $a>L h$
$$
\sup_\zs{S\in \Theta_\zs{\e,L}}\,
\sup_\zs{p\in\cP_\zs{\varsigma}}\,
\P_\zs{S,p}\left(|\wh{\Delta}_\zs{\nu}| > a\right) 
\leq \M_\zs{1,q}\,(\ln n)^{q}\,h^{(1-\alpha)q} +
\frac{\M_\zs{2,q}}{[\iota(a-Lh)^2]^{q/2}}\,,
 $$
where
$\M_\zs{1,q}=2^{q}\,\varrho^{*}_\zs{q}$
and $\M_\zs{2,q}=2^{q}\,\b^{*}_\zs{q}\,\s^{*}_\zs{q}\,r^{*}_\zs{q}$.
 \end{proposition}

\proof
First, we write the estimation error as follows
$$
\wh{\Delta}_\zs{\nu}
 =   B_\zs{\nu}\, + 
\frac{1}{A_\zs{\nu}}\, \zeta_\zs{\nu}\,,
$$
where $\zeta_\zs{\nu} =  \sum^{\nu}_\zs{j=1}\,Q(u_\zs{j})\,y_\zs{j-1}\,\xi_\zs{j}$ and
$$ B_\zs{\nu} =  \frac{1}{A_\zs{\nu}}\,\sum^{\nu}_\zs{j=1}\,Q(u_\zs{j})\,(S(x_\zs{j})-S(z_\zs{0}))\,y^2_\zs{j-1}\,. 
$$  

\noindent Note that $|B_\zs{\nu}| \leq Lh  $ for any $S\in \Theta_\zs{\e,L}$.
 Putting $v=\iota/2 $  we can write

\begin{align}\nonumber
\P_\zs{S,p} \left(|\wh{\Delta}_\zs{\nu}| > a \right)
&=\P_\zs{S,p}\left(|\wh{\Delta}_\zs{\nu}| > a
\,,\,A_\zs{\nu} < v\right)
 +\P_\zs{S,p} \left(|\wh{\Delta}_\zs{\nu}| > a
\,,\,A_\zs{\nu} \geq v\right)\\[2mm] \nonumber
&\leq \P_\zs{S,p}\left(A_\zs{\nu}< v \right)+ \P_\zs{S,p}\left( Lh + \frac{|\zeta_\zs{\nu}|}{A_\zs{\nu}} > a, 
A_\zs{\nu}\geq v \right)\\[2mm] \label{sec:Pr.2}
&\leq \P_\zs{S,p}\left(A_\zs{\nu}< v \right)+ 
\P_\zs{S,p}\left(\frac{|\zeta_\zs{\nu}|}{A_\zs{\nu}} > a-Lh , A_\zs{\nu}\geq v \right)\,.
\end{align}

\noindent 
Now, for any $\bbr\to\bbr$ function $f$ and  numbers $0\le k\le m-1$ 
we set
\begin{equation}\label{A.8}
\varrho_\zs{k,m}(f) = \frac{1}{nh}\sum_\zs{j=k+1}^{m}\,
f(u_\zs{j})\,y^{2}_\zs{j-1}-
\frac{1}{\gamma(S)}\,\frac{1}{nh}\sum_\zs{j=k+1}^{m}\,f(u_\zs{j})\,.
\end{equation}

\noindent Using this function we can estimate 
the first term on the left-hand side of \eqref{sec:Pr.2} as 
\begin{align*}
\P_\zs{S,p}(A_\zs{\nu}< v)
&= \P_\zs{S,p} \left( \varrho_\zs{k_\zs{*},\nu}(Q) + 
\frac{1}{nh\,\gamma(S)} \sum^{\nu}_\zs{j=k_\zs{*}+1}\,Q(u_\zs{j}) < \frac{v}{nh} 
\right)\\
&=
\P_\zs{S,p} \left( \varrho_\zs{k_\zs{*},\nu}(Q) 
< - \frac{\iota}{2 nh} \right)\\
&\leq \P_\zs{S,p} \left( |\varrho_\zs{k_\zs{*},\nu}(Q)| > 
\wt{\epsilon}/2\right) 
\leq
\,
 \left(
\frac{2}{\wt{\epsilon}}
\right)^{q}\,\E_\zs{S,p}|\varrho_\zs{k_\zs{*},\nu}(Q)|^{q} \,.
\end{align*}
Therefore, using here Lemma~\ref{Le.A.3} we get
$$
\P_\zs{S,p}(A_\zs{\nu}< v)\le 2^{q}\,\varrho^{*}_\zs{q}\,
(\ln n)^{q}\,h^{(1-\alpha)q}\,,
$$
where the coefficient  $\varrho^{*}_\zs{q}$ is defined in \eqref{sec:Pr.1}.
The last term on the right-hand side of \eqref{sec:Pr.2} can be estimated as 
\begin{align*}
\P_\zs{S,p}\left(\frac{1}{A_\zs{\nu}}|\zeta_\zs{\nu}| > a-Lh\,,
\,A_\zs{\nu}\geq v \right)
&\leq \P_\zs{S,p}\left(| \zeta_\zs{\nu}| > v(a-Lh) \right)\\
&\leq \frac{1}{v^{q}(a-Lh)^{q}}\,\E_\zs{S,p}|\zeta_\zs{\nu}|^{q}\,. 
\end{align*}
 
\noindent Now in view of the Burkh\"{o}lder inequality, it comes
\begin{align*}
\E_\zs{S,p}|\zeta_\zs{\nu}|^{q} 
&= \E_\zs{S,p}  \left(\sum^{\nu}_\zs{j=1}\,Q(u_\zs{j})\,y_\zs{j-1}\,\xi_\zs{j}\right)^{q} 
\leq \b^{*}_\zs{q}\,\E_\zs{S,p}\,
 \left(\sum^{\nu}_\zs{j=1}\,
Q^2(u_\zs{j})\,y^2_\zs{j-1}\,\xi^2_\zs{j}\right)^{q/2}\\[2mm]
&=\b^{*}_\zs{q}\,\E_\zs{S,p}\,
 \left(\sum^{k_\zs{*}+\iota}_\zs{j=k_\zs{*}+1}
\,y^2_\zs{j-1}\,\xi^2_\zs{j}\right)^{q/2}
\,,
\end{align*}
and after applying the H\"{o}lder inequality, we obtain   
$$
\E_\zs{S,p}|\zeta_\zs{\nu}|^{q} 
 \leq \b^{*}_\zs{q}\,\iota^{q/2-1}\,
\sum^{k_\zs{*}+\iota}_\zs{j=k_\zs{*}+1}\,
\E_\zs{S,p}\, y^{q}_\zs{j-1}\,\xi^{q}_\zs{j} 
\leq \b^{*}_\zs{q}\,\s^{*}_\zs{q}\,r^{*}_\zs{q} \iota^{q/2}\,.
$$
Therefore, 
$$
\P_\zs{S,p}\left(\frac{1}{A_\zs{\nu}}|\zeta_\zs{\nu}| > a-Lh , 
A_\zs{\nu}\geq v \right) \leq 
\frac{2^{q}\,\b^{*}_\zs{q}\,\s^{*}_\zs{q}\,r^{*}_\zs{q}}{\iota^{q/2}(a-Lh)^{q} }\,.
$$
Hence Proposition \ref{Le.S.1}\endproof

\begin{proposition}\label{Le.S.2}
Let the bandwidth $h$ be defined by the conditions 
\eqref{sec:Sp.9}--\eqref{sec:Sp.10}. Then, for all $m\geq 1$ and  $0<\delta<1$ 
$$ 
\lim_\zs{n \to +\infty} 
\sup_\zs{\underline{\alpha}\leq \alpha\leq \alpha^*}\, h^{-m}\,
\sup_\zs{S\in \Theta_\zs{\e,L}}\,
\sup_\zs{p\in\cP_\zs{\varsigma}}\,
\P_\zs{S,p} \left(|\wh{\Delta}_\zs{\nu}| >\delta \wt{\epsilon}\right) = 0\,. 
$$
\end{proposition}

\proof
By applying Proposition \ref{Le.S.1} for $a = \delta \wt{\epsilon}$, we obtain 
that for sufficiently large $n\ge 1$, for which $\delta \wt{\epsilon}>2Lh$,
and for any $m\ge 1$ and $q> m/(1-\overline{\alpha})$
\begin{align*}
 h^{-m}\P_\zs{S,p} (|\wh{\Delta}_\zs{\nu}| > \delta \wt{\epsilon})
& \leq h^{-m}\M_\zs{1,q}\,(\ln n)^{q}h^{(1-\alpha)q} + 
 \frac{\M_\zs{2,q} h^{-m}}
{\left(\iota(\delta\wt{\epsilon} -Lh)^{2}\right)^{q/2}}\\[2mm]
& \leq \M_\zs{1,q}\,(\ln n)^{q}\, \kappa_\zs{n}^{\frac{(1-\overline{\alpha})q-m}{2\overline{\beta}+1}} + 
\frac{2^{q}\,\M_\zs{2,q}}{\delta^{q}}\,\frac{ h^{-(m+q/2)}}
{n^{q/2}\wt{\epsilon}^{3q/2}_\zs{n}}\\[2mm]
& \leq  \M_\zs{1,q}\,(\ln n)^{q} 
\kappa_\zs{n}^{\frac{(1-\overline{\alpha})q-m}{2\overline{\beta}+1}} + 
\frac{2^{q}\,\M_\zs{2,q}(\ln n)^{3q/2}}{\delta^{q}(\kappa_\zs{n}n)^{q/2}}\,
\kappa_\zs{n}^{\frac{q(1-\overline{\alpha}/2)-m}{2\overline{\beta}+1}}\,.
\end{align*}
Taking into account here the conditions \eqref{sec:Sp.10} we come 
to Proposition \ref{Le.S.2}.
\endproof 

\section{Properties of stopping time $\tau_\zs{H}$}\label{sec: Stp}

First  we need to study  
 some asymptotic properties of the term \eqref{sec:Sp.3-0}.

\begin{proposition}\label{Pr.S.3}
Assume that the threshold $H$ is chosen in the form \eqref{sec:Sp.8}
and the bandwidth $h$ satisfies the conditions \eqref{sec:Sp.9} - 
\eqref{sec:Sp.10}. Then for any $m\ge 1$
$$
\limsup_\zs{n\to\infty}\,\sup_\zs{\underline{\alpha}\le \alpha\le \overline{\alpha}}\,h^{-m}
\sup_\zs{S\in \Theta_\zs{\e,L}}\sup_\zs{p\in\cP_\zs{\varsigma}}\,
\P_\zs{S,p} (A_\zs{\nu,n}<H)<\infty\,.
$$
\end{proposition}

\proof
Using the definition of $H$ in \eqref{sec:Sp.8}
we obtain
\begin{align*}
 \P_\zs{S,p} (A_\zs{\nu,n}<H)
&= \P_\zs{S,p}\left(\frac{1}{nh}\sum_\zs{j=\nu+1}^{n} Q(u_\zs{j})\,y^2_\zs{j-1} < \frac{H}{nh}\right)  \\[2mm]
&= \P_\zs{S,p} \left(\varrho_\zs{\nu,k^{*}}(Q)+\frac{1}{\gamma(S)}\,\sum_\zs{j=\nu+1}^{k^*} Q(u_\zs{j})\,\Delta u_\zs{j} < \phi \right)\,.
\end{align*}
Note that 
$$
 \sum_\zs{j=\nu+1}^{k^*} Q(u_\zs{j})\,\Delta u_\zs{j}  
=  \frac{k^{*}-k_\zs{*}-\iota}{nh} \geq  2- \frac{\iota+2}{nh}\,. 
 $$
\noindent 
Taking into account that $\e^2\leq \gamma(S)\leq 1$, we obtain 
\begin{equation} \label{sec:Stp.1}
\left|\frac{1}{\gamma(\wh{S}_\zs{\nu})}- \frac{1}{\gamma(S)}\right|
\leq\frac{2}{\e^4} |\wh{\Delta}_\zs{\nu}|\,. 
\end{equation}
This yields
\begin{align*} 
\P_\zs{S,p}\left(A_\zs{\nu,n}<H\right) 
&\le \P_\zs{S,p}\left(
\sum^{k^{*}}_\zs{j=\nu+1}\,Q(u_\zs{j}) y^{2}_\zs{j-1}<H\,,\,
|\wh{\Delta}_\zs{\nu}| \leq \delta \wt{\epsilon} \right) 
+ \P_\zs{S,p} \left(|\wh{\Delta}_\zs{\nu}| > \delta \wt{\epsilon}\right) 
\\[2mm]
& \leq \P_\zs{S,p} \left(\varrho_\zs{\nu,k^{*}}(Q) <  -2\wt{\epsilon}_n + \frac{\wt{\epsilon}_n}{\e^2} + \frac{4}{\e^4}\delta \wt{\epsilon} + \frac{3}{\e^2 nh} \right)\\[2mm] 
& + 
\P_\zs{S,p} \left(
|\wh{\Delta}_\zs{\nu}| > \delta \wt{\epsilon}
\right)\,.
\end{align*}
Therefore for $\delta<\e^{4}/8$ and sufficiently large $n\ge 1$ we obtain that
$$
\P_\zs{S,p}\left(A_\zs{\nu,n}<H \right) \le 
\P_\zs{S,p} \left(|\varrho_\zs{\nu,k^{*}}(Q)|> \wt{\epsilon}_n/2\right)+
\P_\zs{S,p} \left(
|\wh{\Delta}_\zs{\nu}| > \delta \wt{\epsilon}
\right)\,.
$$
 Lemma~\ref{Le.A.3} and Proposition~\ref{Le.S.2}
imply Proposition \ref{Pr.S.3}.
\endproof

\noindent Now for any weighted sequence $(w_\zs{j})_\zs{j\ge 1}$ we set

\begin{equation}\label{sec:Stp.2}
Z_\zs{n}
= \sum_\zs{j=\tau_\zs{H}-1}^{n}\,Q(u_\zs{j})\,w_\zs{j}\, y^2_\zs{j-1} + 
(1-\varkappa_\zs{H})\,
w_\zs{\tau_\zs{H}}\,Q(u_\zs{\tau_\zs{H}})\,
y^2_\zs{\tau_\zs{H}-1}\,.
\end{equation}

\begin{proposition}\label{Pr.A.4}
Assume that the threshold $H$ is chosen in the form \eqref{sec:Sp.8}
and the bandwidth $h$ satisfies the conditions \eqref{sec:Sp.9} - 
\eqref{sec:Sp.10}. Moreover, let $(w_\zs{j})_\zs{j\ge 1}$ be a
sequence bounded by a constant 
$w^{*}$, i.e. $\sup_\zs{j\ge 1}|w_\zs{j}|\le w^{*}$.
Then 
\begin{equation}\label{sec:Stp.3}
\limsup_\zs{n\to \infty}
\sup_\zs{\underline{\alpha} \leq\alpha\leq\overline{\alpha}}\,h^{-\alpha}\,
\sup_\zs{S\in \Theta_\zs{\e,L}}
\,\sup_\zs{p\in\cP_\zs{\varsigma}}\,\frac{\E_\zs{S,p}|Z_\zs{n}|}{nh} =0\,.
\end{equation}
\end{proposition}
\proof
It is clear that $Z_\zs{n}=0$ if $ A_\zs{\nu,n} < H$, and 
on the set $\{ A_\zs{\nu,n} \geq H\}$ this  term can be estimated as
$$
|Z_\zs{n}| 
\leq   w^{*}
\left(
\sum_\zs{j=\nu+1}^{n}\,Q(u_\zs{j})\, y^2_\zs{j-1} - H
\right)=
 w^{*} (A_\zs{\nu,n} -H ) 
\,,
$$
\noindent 
i.e. $|Z_\zs{n}| 
\leq   w^{*}(A_\zs{\nu,n} -H )_\zs{+}$, where $(x)_\zs{+}=\max(0,x)$.
Therefore,
\begin{align*}
\frac{|Z_\zs{n}|}{nh}&
\leq  w^{*}
\left|\frac{\sum_\zs{j=\nu+1}^{n}\,Q(u_\zs{j})\, y^2_\zs{j-1}}{nh}  
-\phi \right|
 \\[2mm]
&\leq |\varrho_\zs{\nu,n}(Q)| + 
\left|
\frac{\sum_\zs{j=\nu+1}^{n}\,Q(u_\zs{j})}{\gamma(S)nh}\, -\phi
\right|
\,.
\end{align*}

\noindent Taking into account that 
$ \sum_\zs{j=\nu+1}^{n}\,Q(u_\zs{j}) = 
k^{*}-k_\zs{*}-\iota \leq 2nh
$
we obtain 

\begin{align*}
\frac{\E_\zs{S,p}|Z_\zs{n}|}{nh}
& \leq  \E_\zs{S,p}\,|\varrho_\zs{\nu,k^*}(Q)|
+ 2\left|\frac{1}{\gamma(S)}- \frac{1}{\gamma(\wh{S}_\zs{\nu})}\right| + \frac{2}{\e^2}\,\wt{\epsilon}\\
&\leq  \E_\zs{S,p}\,|\varrho_\zs{\nu,k^*}(Q)| + \frac{4}{\e^4}\E_\zs{S,p}\left|
\hat{\Delta}_\zs{\nu}(z_\zs{0})\right| + \frac{2}{\e^2}\,
\wt{\epsilon}\,.
\end{align*}
Moreover, note that
\begin{align*}
\E_\zs{S,p}\left|\hat{\Delta}_\zs{\nu}(z_\zs{0})\right|
&= \E_\zs{S,p}\left|\hat{\Delta}_\zs{\nu}(z_\zs{0})\right|\,
\Chi_\zs{\{\left| \hat{\Delta}_\zs{\nu}(z_\zs{0})\right|\leq \wt{\epsilon}\}}\\[2mm]
&+ 
\E_\zs{S,p}\left|\hat{\Delta}_\zs{\nu}(z_\zs{0})\right|\,
\Chi_\zs{\{\left| \hat{\Delta}_\zs{\nu}(z_\zs{0})\right|> \wt{\epsilon}\}}
\\[2mm]
&\leq \wt{\epsilon} + 2 
\P_\zs{S,p}\left(\left|\hat{\Delta}_\zs{n_0}(z_\zs{0})\right| >  \wt{\epsilon}\right)
\,.
\end{align*}
Therefore,Lemma~\ref{Le.A.3} and Proposition \ref{Le.S.1} imply immediately
\eqref{sec:Stp.3}.  \endproof

\section{Proofs}\label{sec:Prf}

\subsection{Proof of Theorem~\ref{Th.sec:Mr.1}}

First, similarly to the proof of Theorem 2.1 from \cite{ArPe08} we choose the 
corresponding parametric functional
family $S_\zs{u,\delta}(\cdot)$ in the following form 
\begin{equation}\label{sec:Prf.1}
S_\zs{u,\delta}(x)= \frac{u}{\varphi_n}\,
V_\zs{\delta}\left(\frac{x-z_\zs{0}}{h}\right)\,,
\end{equation}
with the function $V_\zs{\delta}$ defined as
$$
V_\zs{\delta}(x)=\delta^{-1}\int^\infty_\zs{-\infty}\wt{Q}_\delta(u) g\left(
\frac{u-x}{\delta}\right) \d u\,,
$$
 where $\wt{Q}_\delta(u)={\bf 1}_\zs{\{|u|\le 1-2\delta\}}+
2{\bf 1}_\zs{\{1-2\delta\le |u|\le 1-\delta\}}$ with $0<\delta<1/4$ and $g$ is some even nonnegative infinitely differentiable
function such that $g(z)=0$ for $|z|\ge 1$ and $\int^1_\zs{-1}\,g(z)\ \d z=1$. 
One can show (see \cite{GaPe06}) that for any $b>0$ and $0<\delta<1/4$ 
there exists $n_\zs{*}=n_\zs{*}(b,L,\delta)>0$ such that for all
$|u|\le b$ and $n\ge n_\zs{*}$
$$
S_\zs{u,\delta} \in 
 \cU^{(\beta)}_\zs{n}(\e,L,\epsilon^{*}_\zs{n})\,.
$$
Therefore, in this case for any $n\ge n_\zs{*}$
\begin{align*}
\varphi_n\,
\sup_\zs{S\in  \cU^{(\beta)}_\zs{n}(\e,L,\epsilon^{*}_\zs{n})} \gamma^{-1/2}(S)\,\cR_\zs{n}(\wt{S}_\zs{n},S)
&\ge 
\sup_\zs{S\in  \cU^{(\beta)}_\zs{n}(\e,L,\epsilon^{*}_\zs{n})} \gamma^{-1/2}(S)
\E_\zs{S,p_\zs{0}}\,\psi_n(\wt{S}_n,S)\\[2mm]
&
\ge\,
\gamma_\zs{*}(n,b)\,
\frac{1}{2 b}\int^{b}_\zs{-b}\,
\E_\zs{S_\zs{u,\delta},p_\zs{0}}\psi_\zs{n}(\wt{S}_\zs{n},S_\zs{u,\delta})\d u \,,
\end{align*}
where $\gamma_\zs{*}(n,b)=\inf_\zs{|u|\le b}\,\gamma^{-1/2}(S_\zs{u,\delta})$.
The definitions \eqref{sec:Sp.6} and \eqref{sec:Prf.1} imply
that for any $b>0$
$$
\lim_\zs{n\to\infty}\sup_\zs{|u|\le b}\,
|\gamma(S_\zs{u,\delta})-1|\,=\,0\,.
$$
Therefore, by the same way as in the proof of Theorem 2.1 from \cite{ArPe08} we obtain that for any 
$b>0$ and $0<\delta<1/4$
\begin{equation}\label{sec:Prf.2}
\underline{\lim}_\zs{n\to\infty} \,\inf_\zs{\wt{S}}\,
\sup_\zs{S\in  \cU^{(\beta)}_\zs{n}(\e,L,\epsilon^{*}_\zs{n})} \gamma^{-1/2}(S)\,
 \varphi_n\,
\cR_\zs{n}(\wt{S}_\zs{n},S)
\ge I(b,\sigma_\zs{\delta})\,,
\end{equation}
where 
$$
I(b,\sigma_\zs{\delta})  = \frac{\max (1,b-\sqrt{b})}{b}\frac{\sigma_\zs{\delta}}{\sqrt{2\pi}}\int_\zs{-\sqrt{b}}^{\sqrt{b}} e^{-\sigma_\zs{\delta}^2\frac{u^2}{2}} du\,,
$$
with $\sigma^2_\zs{\zs\delta}=\int_\zs{-1}^{1}\,V^2_\delta(u)\,\d u$. It is easy to check that
 $\sigma^2_\zs{\zs\delta} \rightarrow 2 $ as $\delta \rightarrow 0$.
Limiting $b\to\infty$  and $\delta \rightarrow 0$ in \eqref{sec:Prf.2} yield
the inequality \eqref{sec:Mr.3}.
Hence Theorem~\ref{Th.sec:Mr.1}.
\endproof

\subsection{Proof of Theorem~\ref{Th.sec:Mr.2}}
First we set
\begin{equation}\label{sec:Prf.2-0}
\check{\varkappa}_\zs{j}=
\Chi_\zs{\{\tau_\zs{H}\neq j\}}
+
\varkappa_\zs{H}\,\Chi_\zs{\{\tau_\zs{H}=j\}}\,.
\end{equation}
Then taking this into account we can represent the estimate error as
\begin{equation}\label{sec:Prf.3}
 S^{*}_\zs{h}-S(z_\zs{0}) = 
-S(z_\zs{0})\,\Chi_\zs{(A_\zs{\nu,n}< H)} + h^{\beta}\,B_\zs{n}(h)\,
\Chi_\zs{(A_\zs{\nu,n}\geq H)}+ 
\frac{1}{\sqrt{H}}\,\zeta_\zs{n}(h)\,\Chi_\zs{(A_\zs{\nu,n}\geq H)}\,,
\end{equation}
where
$$ 
B_\zs{n}(h) =  \frac{
\sum^{\tau_\zs{H}}_\zs{j=\nu+1}
\,\check{\varkappa}_\zs{j}
\,Q(u_\zs{j})\,
(S(x_\zs{j})-S(z_\zs{0}))\,y^2_\zs{j-1}}{h^{\beta}\,H}
$$ 
and
 $$
 \zeta_\zs{n}(h) =\frac{ \sum^{\tau_\zs{H}}_\zs{j=\nu+1}\,
\check{\varkappa}_\zs{j}
\,Q(u_\zs{j})\,y_\zs{j-1}\,
\xi_\zs{j}}{\sqrt{H}}\,.
$$ 

\noindent First we study the term $B_\zs{n}(h)$. To this end we introduce
$$
 B^{*}_\zs{n}= \sum^{n}_\zs{j=\nu+1}\,Q(u_\zs{j})\,b^{*}_\zs{j}\,
 y^{2}_\zs{j-1}\,, 
\quad
b^{*}_\zs{j}=\frac{S(x_\zs{j})-S(z_\zs{0})}{h}
\,\Chi_\zs{\{k_\zs{*}\le j\le k^{*}\}}\,.
$$
It is clear that for any $S$ from $\cU^{(\beta)}_\zs{n}(\e,L,\epsilon^{*}_\zs{n})$
$$
\sup_\zs{j\ge 1}
|b^{*}_\zs{j}|\le L\,.
$$
Therefore, using Proposition~\ref{Pr.A.4} for the sequence  \eqref{sec:Stp.2}
with $w_\zs{j}=b^{*}_\zs{j}$ we obtain 
\begin{equation}\label{sec:Prf.5-00}
\limsup_\zs{n\to \infty}
\sup_\zs{\underline{\alpha} \leq\alpha\leq\overline{\alpha}}\,h^{-\alpha}\,
\sup_\zs{S\in \cU^{(\beta)}_\zs{n}(\e,L,\epsilon^{*}_\zs{n})}
\sup_\zs{p\in\cP_\zs{\varsigma}}\,
\frac{\E_\zs{S,p}|\sum^{\tau_\zs{H}}_\zs{j=\nu+1}
\,\check{\varkappa}_\zs{j}
\,Q(u_\zs{j})\,
b^{*}_\zs{j}\,y^{2}_\zs{j-1}-B^{*}_\zs{n}|}{nh}=0\,.
\end{equation}

\noindent Moreover, putting
$$
f_\zs{1}(u)=Q(u)\,
\frac{S(z_\zs{0}+hu)-S(z_\zs{0})}{h}
\,,
$$
we obtain
\begin{equation}\label{sec:Prf.5-0}
\frac{B^{*}_\zs{n}}{nh}=
\frac{1}{\gamma(S)}\sum^{n}_\zs{k=1}\,f_\zs{1}(u_\zs{j})\Delta u_\zs{j}
-
\frac{1}{\gamma(S)}\sum^{\nu}_\zs{k=1}\,f_\zs{1}(u_\zs{j})\Delta u_\zs{j}
+
\varrho_\zs{\nu,n}(f_\zs{1})
\,.
\end{equation}
Using the definition of $\Omega_\zs{h}(z_\zs{0},S)$ in \eqref{sec:Sp.1-0}
we can represent the first term  as
$$
\sum^{n}_\zs{k=1}\,f_\zs{1}(u_\zs{j})\Delta u_\zs{j}=
\frac{1}{h}\Omega_\zs{h}(z_\zs{0},S)
-\int^{1}_\zs{u_\zs{k^{*}}}\,f_\zs{1}(u)\d u
+\sum^{k^{*}}_\zs{k=k_\zs{*}}\,\int^{u_\zs{j}}_\zs{u_\zs{j-1}}
\left(
f_\zs{1}(u_\zs{j})
-
f_\zs{1}(u)
\right)\,\d u\,.
$$
Note now that
for any $S$ from $\cU^{(\beta)}_\zs{n}(\e,L,\epsilon^{*}_\zs{n})$
$$
\sup_\zs{-1\le u\le 1}
|f_\zs{1}(u)|\le L
\quad\mbox{and}\quad
\sup_\zs{-1\le u\le 1}
|\dot{f}_\zs{1}(u)|\le L\,,
$$
and, therefore,
$$
\left|
\int^{1}_\zs{u_\zs{k^{*}}}\,f_\zs{1}(u)\d u
\right|\le \frac{L}{nh}
\quad\mbox{and}\quad
\left|
\sum^{k^{*}}_\zs{k=k_\zs{*}}\,\int^{u_\zs{j}}_\zs{u_\zs{j-1}}
\left(
f_\zs{1}(u_\zs{j})
-
f_\zs{1}(u)
\right)\,\d u
\right|\le 
\frac{2L}{nh}\,.
$$
The last bounds imply immediately
$$
\limsup_\zs{n\to \infty}
\sup_\zs{\underline{\alpha} \leq\alpha\leq\overline{\alpha}}\,h^{-\alpha}\,
\sup_\zs{S\in \cU^{(\beta)}_\zs{n}(\e,L,\epsilon^{*}_\zs{n})}
\,\left|
\sum^{n}_\zs{k=1}\,f_\zs{1}(u_\zs{j})\Delta u_\zs{j}
\right|
=0\,.
$$
Taking into account Lemma~\ref{Le.A.3} in \eqref{sec:Prf.5-0}
we get
$$
\limsup_\zs{n\to \infty}
\sup_\zs{\underline{\alpha} \leq\alpha\leq\overline{\alpha}}\,h^{-\alpha}\,
\sup_\zs{S\in \cU^{(\beta)}_\zs{n}(\e,L,\epsilon^{*}_\zs{n})}\,
\sup_\zs{p\in\cP_\zs{\varsigma}}\,\frac{\E_\zs{S,p}|B^{*}_\zs{n}|}{nh}=0\,.
$$
Therefore, in view of \eqref{sec:Prf.5-00} 
$$
\limsup_\zs{n\to \infty}
\sup_\zs{\underline{\alpha} \leq\alpha\leq\overline{\alpha}}\,
\sup_\zs{S\in \cU^{(\beta)}_\zs{n}(\e,L,\epsilon^{*}_\zs{n})}\,
\sup_\zs{p\in\cP_\zs{\varsigma}}\,\E_\zs{S,p}\,|B_\zs{n}(h)\,=0\,.
$$

\noindent 
To study the last term in \eqref{sec:Prf.2-0} note that 
 the definition of the stopping time in \eqref{sec:Sp.4} implies 
\begin{equation}\label{sec:Up.3}
\sup_\zs{n\ge \nu+1}\,
\sup_\zs{h_\zs{*}\le h\le h^*}\,
\sup_\zs{S\in\Theta_\zs{\e,L}}\,
\sup_\zs{p\in\cP_\zs{\varsigma}}
\E_\zs{S,p}\,|\zeta_\zs{n}(h)|^{2} \le 1\,.
\end{equation}
Therefore, in view of Lemma~\ref{Le.A.5} we obtain
$$
\lim_\zs{n\to\infty}\,
\sup_\zs{S\in\Theta_\zs{\e,L}}\,
\sup_\zs{p\in\cP_\zs{\varsigma}}
\left|
\E_\zs{S,p}\,|\zeta_\zs{n}(h)|\Chi_\zs{(A_\zs{\nu,n}\geq H)}
-
\E |\zeta_\zs{\infty}|
\right|=0\,,
$$
where $\zeta_\zs{\infty}\sim \cN(0,1)$, i.e. $\E |\zeta_\zs{\infty}|=\sqrt{2/\pi}$. 
Moreover, in view of
Proposition~\ref{Le.S.2} and the bound \eqref{sec:Stp.1} 
we get 
$$
\lim_\zs{n\to\infty}
\sup_\zs{S\in\Theta_\zs{\e,L}}\,
\sup_\zs{p\in\cP_\zs{\varsigma}}
\E_\zs{S,p}
\left|
\frac{\varphi^{2}_\zs{n}}{\gamma(S)\,H}
-\frac{1}{2}
\right|
=0\,.
$$
From this and Proposition~\ref{Pr.S.3} it follows Theorem~\ref{Th.sec:Mr.2}. \endproof

\medskip

\subsection{Proof of Theorem~\ref{Th.sec:Mr.1_1}}

First, note that the representation \eqref{sec:Prf.3}
implies
\begin{align}
\E_\zs{S,p}|S^{*}_\zs{h}-S(z_\zs{0})|
& \leq \P_\zs{S,p}(A_\zs{\nu,n}< H) + h^{\beta}\,
\E_\zs{S,p}|B_\zs{n}(h)|\,\Chi_\zs{(A_\zs{\nu,n}\geq H)} 
\nonumber\\[2mm]\label{sec:Prf.6}
& + \E_\zs{S,p}\left(
\frac{|\zeta_\zs{n}(h)|}{\sqrt{H}}\,\,\Chi_\zs{(A_\zs{\nu,n}\geq H)}
\right)\,.
\end{align}
Let us show, that
\begin{equation}\label{sec:Prf.7}
\limsup_\zs{n\to\infty}\,\sup_\zs{\underline{\alpha}\leq \alpha \leq \overline{\alpha}}\,
\sup_\zs{S\in\cH^{(\beta)}(\e,L,L^{*})}\,\sup_\zs{p\in \cP}\,
\E_\zs{S,p}\,|B_\zs{n}(h)| < \infty\,.
\end{equation}
Indeed, setting 
\begin{equation}\label{sec:Prf.3-1}
\vartheta_\zs{j} = \frac{S(x_\zs{j})-S(z_\zs{0})}{h} - 
\dot{S}(z_\zs{0})\,u_\zs{j}\,,
\end{equation}
we can represent $B_\zs{n}(h)$ as
\begin{equation}\label{sec:Prf.3-2}
B_\zs{n}(h) 
= \frac{h^{-\alpha}}{H}\dot{S}(z_\zs{0}) \wt{B}_\zs{n} +
\frac{h^{-\alpha}}{H} \wh{B}_\zs{n}\,,
\end{equation}
where
$$
\wt{B}_\zs{n}= \displaystyle\sum^{\tau_\zs{H}}_\zs{j=\nu+1}
\,\check{\varkappa}_\zs{j}
\,Q(u_\zs{j})\,u_\zs{j}\, y^2_\zs{j-1}
\quad\mbox{and}\quad
\wh{B}_\zs{n}= \displaystyle\sum^{\tau_\zs{H}}_\zs{j=\nu+1}\,
\check{\varkappa}_\zs{j}
\,
Q(u_\zs{j})\,\vartheta_\zs{j} \, y^2_\zs{j-1}\,. 
$$

\noindent
Using now Proposition~\ref{Pr.A.4} for the sequence  \eqref{sec:Stp.2}
with $w_\zs{j}=u_\zs{j}$ we obtain that
$$
\limsup_\zs{n\to \infty}
\sup_\zs{\underline{\alpha} \leq\alpha\leq\overline{\alpha}}\,h^{-\alpha}\,
\sup_\zs{S\in \Theta_\zs{\e,L}}\,
\sup_\zs{p\in\cP_\zs{\varsigma}}\,\frac{\E_\zs{S,p}|\wt{V}_\zs{n}-\wt{B}_\zs{n}|}{nh} =0\,,
$$
where 
$$
 \wt{V}_\zs{n}= \sum^{n}_\zs{j=\nu+1}\,Q(u_\zs{j})\,u_\zs{j}\, y^2_\zs{j-1}\,.
$$
We can represent this term as 
\begin{align*}
\wt{V}_\zs{n}
&= \varrho_\zs{\nu,n}(Q_\zs{1})+ \frac{1}{\gamma (S)}\,
\sum^{n}_\zs{j=\nu+1}\,Q_\zs{1}(u_\zs{j})\,\Delta u_\zs{j}\nonumber \\
&= \varrho_\zs{\nu,n}(Q_\zs{1})+
\frac{1}{\gamma (S)}\,\left(\sum^{k^{*}}_\zs{j=k_\zs{*}+1}\,u_\zs{j}\,
\Delta u_\zs{j} - \sum^{k_\zs{*}+
\iota}_\zs{j=k_\zs{*}+1}\,u_\zs{j}\,\Delta u_\zs{j}\right)
\,,
\end{align*}
where $Q_\zs{1}(u) = Q(u)\,u$. Moreover, taking into account here, that
\begin{align*}
\left|\sum^{k^{*}}_\zs{j=k_\zs{*}+1}\,u_\zs{j}\,\Delta u_\zs{j}\right|
&=
\left|\int^{u_\zs{k^{*}}}_\zs{u_\zs{k_\zs{*}}}\,u\,\d u
+
\sum^{k^{*}}_\zs{j=k_\zs{*}+1}\,\int^{u_\zs{j}}_\zs{u_\zs{j-1}}\,
\left(u_\zs{j}-u\right)\,\d u
\right|
 \leq \frac{4}{nh}
\end{align*}
and
\begin{align*}
\left|
\sum^{k_\zs{*}+\iota}_\zs{j=k_\zs{*}+1}\,u_\zs{j}\,\Delta u_\zs{j}
\right| \leq \wt{\epsilon}+\frac{1}{nh}
\,,
\end{align*}
we obtain
$$ 
\frac{1}{nh}\E_\zs{S,p}\,|\wt{V}_\zs{n}| \leq 
\E_\zs{S,p}\,|\varrho_\zs{\nu,n}(Q_\zs{1})|
 + \frac{1}{\e^2}\left(\frac{5}{nh} + \wt{\epsilon}\right)
\,.
$$
Now Lemma~\ref{Le.A.3} yields
$$
\limsup_\zs{n\to \infty}\sup_\zs{\underline{\alpha} 
\leq\alpha\leq\overline{\alpha}}\,h^{-\alpha}\,
\,
\sup_\zs{S\in\Theta_\zs{\e,L}}\,
\sup_\zs{p\in\cP_\zs{\varsigma}}\frac{\E_\zs{S,p}|\wt{V}_\zs{n}|}{nh} =0 
$$
and, therefore,
\begin{equation}\label{sec:Prf.5}
\limsup_\zs{n\to \infty}
\sup_\zs{\underline{\alpha} \leq\alpha\leq\overline{\alpha}}\,h^{-\alpha}\,
\sup_\zs{S\in \Theta_\zs{\e,L}}
\sup_\zs{p\in\cP_\zs{\varsigma}}\,\frac{\E_\zs{S,p}|\wt{B}_\zs{n}|}{nh} =0
\,.
\end{equation}
\noindent To estimate the last term in \eqref{sec:Prf.3-2}
note, that for any function $S$ from $\cH^{(\beta)}(\e,L,L^{*})$
and for $k_\zs{*}\le j\le k^{*}$ the coefficients \eqref{sec:Prf.3-1}
can be estimated as
$$
 |\vartheta_\zs{j}|=\left| \int^{u_\zs{j}}_\zs{0}
\left(\dot{S}(z_\zs{0}+hu) 
- \dot{S}(z_\zs{0})\right)\d u\right| \, 
\le L|u_\zs{j}|\, h^{\alpha} 
\le L^{*}\,h^{\alpha}\,.
$$
Therefore, 
$$
\limsup_\zs{n\to\infty}\,\sup_\zs{\underline{\alpha}\le \alpha \le \overline{\alpha}}\,h^{-\alpha}\,
\sup_\zs{S\in\cH^{(\beta)}(\e,L,L^{*})}\,\sup_\zs{p\in \cP}
\, \frac{1}{n h}
\E_\zs{S,p}\,|\wh{B}_\zs{n}|\,<\infty\,.
$$
Now the property \eqref{sec:Prf.5} implies the  inequality \eqref{sec:Prf.7}. 
Therefore, using  Proposition \ref{Pr.S.3} and the inequality \eqref{sec:Up.3} in
\eqref{sec:Prf.6},  we come to Theorem \ref{Th.sec:Mr.2}.
\endproof

\medskip

\subsection{Proof of Theorem~\ref{Th.6.2}}

First of all, note that 
the coefficient $\phi$ defined in  \eqref{sec:Sp.8}
will be more than one for sufficient large $n$ for which $\wt{\epsilon}\le 1/2$.
So, using the representation \eqref{sec:Prf.3}, we get 
for any $1\le j\le m$
\begin{equation}\label{sec:Prf.8}
|\check{S}_\zs{j}-S(z_\zs{0})|\le  
\Chi_\zs{\Gamma_\zs{j}} + 
(\check{h}_\zs{j})^\beta\, |B_\zs{n}(\check{h}_\zs{j})|+ 
\frac{1}{\sqrt{n\check{h}_\zs{j}}}
\,|\wt{\zeta}_\zs{n}(\check{h}_\zs{j})|\,,
\end{equation}
where $\Gamma_\zs{j}=\{A_\zs{\nu,n}(\check{h}_\zs{j})< \check{H}_\zs{j}\}$,
$\wt{\zeta}_\zs{n}(\check{h}_\zs{j})=
\wt{\zeta}_\zs{n}(\check{h}_\zs{j})\,\Chi_\zs{(A_\zs{\nu,n}(\check{h}_\zs{j})\geq \check{H}_\zs{j})}$,
the random functions $B_\zs{n}(h)$ and $\zeta_\zs{n}(h)$
are defined in \eqref{sec:Prf.3}. We set   
$$
i_\zs{0}=\left[\frac{m(\beta- \underline{\beta})}{\overline{\beta}-\underline{\beta}}
\right]
\,.
$$
This means that
$$
0\le \beta-\beta_\zs{i_\zs{0}}< 
\frac{\overline{\beta}-\underline{\beta}}{m}
\,.
$$
Therefore, taking into account the definition 
of $m$ in  \eqref{sec:Ad.2}, we obtain that 
 for any fixed integer $-<\infty<l<\infty$
\begin{equation}\label{sec:Prf.9}
\left\{
\begin{array}{cc}
0<\liminf_\zs{n\to\infty}&\dfrac{N(\beta)}{N_\zs{i_\zs{0}+l}}
\le \limsup_\zs{n\to\infty}\dfrac{N(\beta)}{N_\zs{i_\zs{0}+l}}<\infty\,;
\\[6mm]
0<\liminf_\zs{n\to\infty}&\dfrac{h(\beta)}{\check{h}_\zs{i_\zs{0}+l}}
\le \limsup_\zs{n\to\infty}\dfrac{h(\beta)}{\check{h}_\zs{i_\zs{0}+l}}<\infty
\,.
\end{array}
\right.
\end{equation}
These inequalities and Theorem~\ref{Th.sec:Mr.1_1} imply
\begin{equation}\label{sec:Prf.10}
\limsup_\zs{n\to\infty}\,
\sup_\zs{\underline{\beta}\le \beta\le\overline{\beta}}\,N(\beta)\,
\sup_\zs{S\in  \cH^{(\beta)}(\e,L,L^{*})}\,
\E_\zs{S}\,
\varpi(i_\zs{0},z_\zs{0})\,
\,<\infty\,,
\end{equation}
where
$\varpi(i_\zs{0},z_\zs{0})=|\check{S}_\zs{i_\zs{0}-1}-S(z_\zs{0})|
+|\check{S}_\zs{i_\zs{0}}-S(z_\zs{0})|+|\check{S}_\zs{i_\zs{0}+1}-S(z_\zs{0})|$.
Now considering the estimator $\wh{S}_\zs{a,n},$ one has
\begin{equation}\label{sec:Prf.11}
|\wh{S}_\zs{a,n}-S(z_\zs{0})|\le
I_\zs{1}+I_\zs{2}
+\varpi(i_\zs{0},z_\zs{0})\,,
\end{equation}
where
$
I_\zs{1}=|\wh{S}_\zs{a,n}-S(z_\zs{0})|\Chi_\zs{\{\check{k}\ge i_\zs{0}+2\}}
\quad\mbox{and}\quad
I_\zs{2}=|\wh{S}_\zs{a,n}-S(z_\zs{0})|\Chi_\zs{\{\check{k}\le i_\zs{0}-2\}}$.
We start with the term $I_\zs{1}$. We have
$$
|\wh{S}_\zs{a,n}-S(z_\zs{0})|\Chi_\zs{\{\check{k}\ge i_\zs{0}+2\}}\le
|\wh{S}_\zs{a,n}-\check{S}_\zs{i_\zs{0}}|\Chi_\zs{\{\check{k}\ge i_\zs{0}+2\}}
+
|
\check{S}_\zs{i_\zs{0}}
-S(z_\zs{0})|\Chi_\zs{\{\check{k}\ge i_\zs{0}+2\}}\,.
$$
Moreover,
\begin{align*}
|\wh{S}_\zs{a,n}&-\check{S}_\zs{i_\zs{0}}|\,\Chi_\zs{\{\check{k}\ge i_\zs{0}+2\}}\,
\le
\check{\omega}_\zs{\check{k}}\,\Chi_\zs{\{\check{k}\ge i_\zs{0}+2\}}+
\frac{\check{\lambda}}{N_\zs{i_\zs{0}}}\\
&\le
\frac{\check{\lambda}}{N_\zs{\check{k}}}\Chi_\zs{\{\check{k}\ge i_\zs{0}+2\}}
+
\frac{\check{\lambda}}{N_\zs{i_\zs{0}}}
\le
\frac{\check{\lambda}}{N(\beta)}
+
\frac{\check{\lambda}}{N_\zs{i_\zs{0}}}
\,.
\end{align*}
The inequalities
\eqref{sec:Prf.9}--\eqref{sec:Prf.10}
imply immediately
$$
\limsup_\zs{n\to\infty}\,
\sup_\zs{\underline{\beta}\le \beta\le\overline{\beta}}\,N(\beta)\,
\sup_\zs{S\in  \cH^{(\beta)}(\e,L,L^{*})}\,
\E_\zs{S}\,I_\zs{1}\,<\infty\,.
$$
\noindent Now we study the term $I_\zs{2}$. From \eqref{sec:Prf.8}
it follows that
$$
I_\zs{2}\le 
\left(
\Chi_\zs{\Gamma_\zs{\check{k}}} + 
(\check{h}_\zs{\check{k}})^\beta\, |B_\zs{n}(\check{h}_\zs{\check{k}})|+ 
\frac{1}{\sqrt{n\check{h}_\zs{\check{k}}}}
\,|\wt{\zeta}_\zs{n}(\check{h}_\zs{\check{k}})|
\right)
\,\Chi_\zs{\{\check{k}\le i_\zs{0}-2\}}\,.
$$
Therefore,
\begin{equation}\label{sec:Prf.12}
\E_\zs{S}\,I_\zs{2}
\le \,m\,
I^{*}_\zs{1}(S)
+
I^{*}_\zs{2}(S)\,
\sum^{i_\zs{0}-2}_\zs{j=0}\,\check{h}_\zs{j}^{\beta}
+ \Psi_\zs{n}(S)
\,,
\end{equation}
where $I^{*}_\zs{1}(S)= \max_\zs{0\le l\le m}\,
\P_\zs{S}\left(\Gamma_\zs{l}\right)$,
$
I^{*}_\zs{2}(S)=
 \max_\zs{0\le l\le m}
\,\E_\zs{S}|B_\zs{n}(\check{h}_\zs{l})|
$
and
$$
\Psi_\zs{n}(S)=\frac{1}{\sqrt{n}}
\,
\sum^{i_\zs{0}-2}_\zs{j=0}
\frac{1}{\sqrt{\check{h}_\zs{j}}}
\E_\zs{S}\,
|\wt{\zeta}_\zs{n}(\check{h}_\zs{j})|
\,\Chi_\zs{\{\check{k}=j\}}\,.
$$
Note now, that for any $0\le j\le i_\zs{0}$
and for sufficiently large  $n$ (for which $\ln d_\zs{n}\ge m/2$) we get 
$$
N(\beta)\check{h}_\zs{j}^{\beta}\le e^{-\beta^{*}(i_\zs{0}-j)}\,,
\quad
\beta^{*}=\frac{\beta(\overline{\beta}-\underline{\beta})}{(2\beta+1)(2\overline{\beta}+1)}\,,
$$
i.e.
$$
N(\beta)\,\sum^{i_\zs{0}-2}_\zs{j=0}\,\check{h}_\zs{j}^{\beta}\le 
\frac{e^{\beta^{*}}}{e^{\beta^{*}}-1}\,.
$$
So, Proposition~\ref{Pr.S.3} and the inequality \eqref{sec:Prf.7}
yield
$$
\limsup_\zs{n\to\infty}\,
N(\beta)\,
\sup_\zs{S\in  \cH^{(\beta)}(\e,L,L^{*})}\,
\left(m
I^{*}_\zs{1}(S)+I^{*}_\zs{2}(S)\,
\sum^{i_\zs{0}-2}_\zs{j=0}\,\check{h}_\zs{j}^{\beta}
\right)\,<\,\infty\,.
$$
Let us consider now the last term \eqref{sec:Prf.12}. To this end note that
for $0\le j\le i_\zs{0}-2$
$$
\{\check{k}=j\}\subseteq \{\check{\omega}_\zs{j+1}\ge \check{\lambda}/N_\zs{j+1}\}
\subseteq \cup^{j+1}_\zs{l=0}
\{|\check{S}_\zs{l}-S(z_\zs{0})|\ge \check{\lambda}/N_\zs{l}\}\,.
$$ 
Therefore,
$$
\Psi_\zs{n}(S)\le\frac{1}{\sqrt{n}}
\,
\sum^{i_\zs{0}-2}_\zs{j=0}
\frac{1}{\sqrt{\check{h}_\zs{j}}}
\sum^{j+1}_\zs{l=0}\,
\E_\zs{S}\,
|\wt{\zeta}_\zs{n}(\check{h}_\zs{j})|
\,\Chi_\zs{\{
|\check{S}_\zs{l}-S(z_\zs{0})|\ge \check{\lambda}/N_\zs{l}
\}}\,.
$$
Taking into account that $\check{h}_\zs{j+1}/\check{h}_\zs{j}\le e$, we 
can rewrite the last inequality as
\begin{equation}\label{sec:Prf.13}
\Psi_\zs{n}(S)\le \frac{e}{\sqrt{n}}
\,
\sum^{i_\zs{0}-1}_\zs{j=1}
\frac{1}{\sqrt{\check{h}_\zs{j}}}
\sum^{j}_\zs{l=0}\,
\E_\zs{S}\,
|\wt{\zeta}_\zs{n}(\check{h}_\zs{j})|
\,\Chi_\zs{\{
|\check{S}_\zs{l}-S(z_\zs{0})|\ge \check{\lambda}/N_\zs{l}
\}}\,.
\end{equation}

Now, taking into account the inequality \eqref{sec:Up.3}, we get
\begin{align*}
 \E_\zs{S}\,
|\wt{\zeta}_\zs{n}(\check{h}_\zs{j})|
\,\Chi_\zs{\{
|\check{S}_\zs{l}-S(z_\zs{0})|\ge \check{\lambda}/N_\zs{l}
\}}&\le \sqrt{\P_\zs{S}(\Gamma_\zs{l})}+
\E_\zs{S}\,
|\wt{\zeta}_\zs{n}(\check{h}_\zs{j})|
\,\Chi_\zs{\{
\check{h}^{\beta}_\zs{l}|B_\zs{n}(\check{h}_\zs{l})|\ge \check{\lambda}_\zs{1}/ N_\zs{l}
\}}\\[2mm]
&+
\E_\zs{S}\,
|\zeta^{*}|
\,\Chi_\zs{\{
|\zeta^{*}|\ge \sqrt{\ln n} \check{\lambda}_\zs{1}
\}}
\end{align*}
where $\check{\lambda}_\zs{1}=\check{\lambda}/2$ and 
$\zeta^{*}=\max_\zs{1\le j\le m}\,|\wt{\zeta}_\zs{n}(\check{h}_\zs{j})|$.
In view of the H\"older and Chebyshev inequalities
and making use of  the upper bound \eqref{sec:A.10}
 we obtain
$$
\E_\zs{S}\,
|\wt{\zeta}_\zs{n}(\check{h}_\zs{j})|
\,\Chi_\zs{\{
\check{h}^{\beta}_\zs{l}|B_\zs{n}(\check{h}_\zs{l})|\ge \check{\lambda}_\zs{1}/ N_\zs{l}
\}}\le \,
\frac{(\mu^{*}_\zs{4})^{1/4}(I^{*}_\zs{2}(S))^{3/4}}{\check{\lambda}^{3/4}_\zs{1}}\,
\left( N_\zs{l} \check{h}^{\beta}_\zs{l}\right)^{3/4}\,.
$$
where the term $I^{*}_\zs{2}(S)$ is defined in \eqref{sec:Prf.12}. 
Using these bounds in \eqref{sec:Prf.13} we get
\begin{align}\nonumber
\Psi_\zs{n}(S)&\le \frac{e\,m^{2}}{\sqrt{n\,\check{h}_\zs{0}}}\,
\sqrt{I^{*}_\zs{1}(S)}
+
\frac{e(\mu^{*}_\zs{4})^{1/4}(I^{*}_\zs{2}(S))^{3/4}}{\check{\lambda}_\zs{1}^{3/4}}\,
\Upsilon^{*}_\zs{n}\\[2mm]\label{sec:Prf.14}
&+\frac{e\,m^{2}}{\sqrt{n\,\check{h}_\zs{0}}}\,\E_\zs{S}\,
|\zeta^{*}|
\,
\Chi_\zs{\left\{|\zeta^{*}|\ge \check{\lambda}_\zs{1}\sqrt{\ln n} 
\right\}}\,,
\end{align}
where
$$
\Upsilon^{*}_\zs{n}=\frac{1}{\sqrt{n}}
\,
\sum^{i_\zs{0}-1}_\zs{j=1}
\left( N_\zs{j}\check{h}^{\beta}_\zs{j}\right)^{1/4}
\,
\left( N_\zs{j}\check{h}^{\alpha}_\zs{j}\right)^{1/2}
\,
\upsilon^{*}_\zs{j}
\quad\mbox{and}\quad
\upsilon^{*}_\zs{j}=
\sum^{j}_\zs{l=0}\,\left(
\frac{N_\zs{l} \check{h}^{\beta}_\zs{l}}{N_\zs{j}\check{h}^{\beta}_\zs{j}}
\right)^{3/4}
\,.
$$
Let us estimate the term $\Upsilon^{*}_\zs{n}$. To this end note, that for any
$0\le j< i_\zs{0}$ and for sufficiently large  $n$ (for which $\ln d_\zs{n}\ge m/2$)
$$
\left( N_\zs{j}\check{h}^{\beta}_\zs{j}\right)^{1/4}
=\exp\left\{
\frac{\beta_\zs{j}-\beta}{4(2\beta_\zs{j}+1)}\,\ln d_\zs{n}
\right\}
\le e^{-\beta^{*}_\zs{1}(i_\zs{0}-j)}
\quad\mbox{with}\quad
 \beta^{*}_\zs{1}=\frac{\overline{\beta}-\underline{\beta}}{8(2\overline{\beta}+1)}\,,
$$
and
$$
\left(N_\zs{j}\check{h}^{\alpha}_\zs{j}\right)^{1/2}
=\exp\left\{
\frac{\beta_\zs{j}-\alpha}{2(2\beta_\zs{j}+1)}\,\ln d_\zs{n}
\right\}
\,
\le \frac{1}{\sqrt{\check{h}(\beta)}}=
\frac{\sqrt{d_\zs{n}}}{N(\beta)}\,,
$$
where $\check{h}(\beta)$ is defined in \eqref{sec:Ad.1-1}.
Similarly for any $0\le l\le j\le i_\zs{0}$ we get
$$
\left(
\frac{N_\zs{l} \check{h}^{\beta}_\zs{l}}{N_\zs{j}\check{h}^{\beta}_\zs{j}}
\right)^{3/4}
=
\exp\left\{\frac{3(\beta_\zs{l}-\beta_\zs{j})(2\beta+1)}{4(2\beta_\zs{l}+1)(2\beta_\zs{j}+1)}
\ln d_\zs{n}
\right\}
\le e^{-\beta^{*}_\zs{2}(j-l)}
\,,\quad \beta^{*}_\zs{2}=\frac{3(\overline{\beta}-\underline{\beta})}{8(2\beta+1)}\,.
$$ 
This means that the sequence $(\upsilon^{*}_\zs{j})_\zs{j\ge 1}$ is bounded, i.e.
$$
\sup_\zs{j\ge 1}
\upsilon^{*}_\zs{j}\le \frac{e^{\beta^{*}_\zs{2}}}{e^{\beta^{*}_\zs{2}}-1}\,.
$$
Therefore, 
\begin{equation}\label{sec:Prf.15}
\lim_\zs{n\to\infty}\,
N(\beta)\,\Upsilon^{*}_\zs{n}=0\,.
\end{equation}

\noindent
The last term in \eqref{sec:Prf.14} can be estimated through
 Lemma~\ref{sec:A.9}, i.e.
\begin{align*}
\E_\zs{S}\,\zeta^{*}\,
\Chi_\zs{\left\{\zeta^{*}\ge\check{\lambda}_\zs{1}\,\sqrt{\ln n}\right\}}
&\le\,m\, \max_\zs{1\le j\le m}\,
\E_\zs{S}\,|\wt{\zeta}_\zs{n}(\check{h}_\zs{j})|\,
\Chi_\zs{\left\{|\wt{\zeta}_\zs{n}(\check{h}_\zs{j})|\ge\check{\lambda}_\zs{1}\,\sqrt{\ln n}\right\}}\\[2mm]
&\le \,2m\,
\check{\lambda}_\zs{1}\,\sqrt{\ln n}\, 
e^{-\frac{1}{8}\,\check{\lambda}_\zs{1}^2\,\ln n} + 
2\,m \int_\zs{\check{\lambda}_\zs{1}\,\sqrt{\ln n}}^{+\infty}\,e^{-z^2/8}\,\d z\,.
\end{align*}
Therefore, for sufficiently large $n$ (when $\check{\lambda}_\zs{1}\,\sqrt{\ln n}\ge 1$)
we get that
\begin{equation}\label{sec:Prf.16}
\E_\zs{S}\,\zeta^{*}\,
\Chi_\zs{\left\{\zeta^{*}\ge\check{\lambda}_\zs{1}\,\sqrt{\ln n}\right\}}
 \le\,2m \left(\check{\lambda}_\zs{1}\,\sqrt{\ln n} +4\right)\, 
n^{-\check{\lambda}_\zs{1}^2/8}\,.
\end{equation}
Now the definition of the parameter $\check{\lambda}$
in \eqref{sec:Ad.3} yields
$$
\limsup_\zs{n\to\infty}
\frac{N(\beta)\,m^{2}}{\sqrt{n\,\check{h}_\zs{0}}}\,
\sup_\zs{S\in  \cH^{(\beta)}(\e,L,L^{*})}\,
\E_\zs{S}\,
|\zeta^{*}|
\,
\Chi_\zs{\left\{|\zeta^{*}|\ge \check{\lambda}_\zs{1}\sqrt{\ln n} 
\right\}}
=0\,.
$$
Hence Theorem~\ref{Th.6.2}. \endproof
 
\medskip

\section{Numerical examples}\label{sec:Nur}
We illustrate the obtained results by the following simulation which is established using .
\subsection{Nonadaptive estimation}
In this section we illustrate the results obtained in the case of nonadaptive estimation 
The purpose is to estimate, at a given point $z_\zs{0}$, the function $S$ defined 
over $[0;1]$ by 
\begin{equation}\label{sec:Nur.1}
 S(x) = (x-z_\zs{0})\,|x-z_\zs{0}|^{\alpha}
\end{equation}
for $z_\zs{0}=1/\sqrt{2}$ and $\alpha=0,3$. 
Taking into account that for this function
$$
\Omega_\zs{h}(z_\zs{0},S)=0 
$$
we obtain that 
for any 
\begin{equation}\label{sec:Nur.2}
0<\e\le 1-\left(\frac{1}{\sqrt{2}}\right)^{\beta}
\,,\quad L\ge\beta
\quad\mbox{and}\quad
\epsilon^{*}_\zs{n}>0
\end{equation}
the
function
\eqref{sec:Nur.1}
 belongs to the class
$\cU^{(\beta)}_\zs{n}(\e,L,\epsilon^{*}_\zs{n})$ with $\beta=1,3$.

The numerical results approximate the asymptotic risk of
a estimators defined in 
\eqref{sec:Sp.5}
  used due to the calculation of an expectation
(it performs an average for $ M = 30000$ simulations) and the finite number of
observations $n$. Here we calculate for the estimator the quantity 
$$
\R_\zs{n} = \displaystyle\frac{1}{M}\,\sum_\zs{k=1}^{M}\,
\left|S^{*,k}_\zs{h} - 
S(z_\zs{0})
\right|\,.
$$ 

For  the standard Gaussian random variables $(\xi_\zs{j})_\zs{j\ge 1}$ in \eqref{sec:In.1},
and by varying the number of observations $n,$ we obtain different risks listed in 
the following table: we obtain:\\
\begin{center}
\begin{tabular}{|c|c|c|c|c|c|}
  \hline
  $n$ &         1000 &  5000& 10000 & 20000\\
  \hline
 $\R_n$ &  0.034 & 0.021 & 0.017 &0.012  \\
  \hline
  
\end{tabular}
\end{center}

For random variables $(\xi_\zs{j})_\zs{j\ge 1}$ reduced from uniform random variables on $[-1,1],$ we obtain :

\begin{center}
\begin{tabular}{|c|c|c|c|c|}
  \hline
  $n$ &          1000 &  5000& 10000 & 20000 \\
  \hline
 $\R_n$ &  0.038 & 0.022 & 0.018 &0.014 \\
  \hline
  
\end{tabular}
\end{center}

For random variables $(\xi_\zs{j})_\zs{j\ge 1}$ centered and reduced from exponential random variables with parameter $1,$ we obtain :

\begin{center}
\begin{tabular}{|c|c|c|c|c|}
  \hline
  $n$ &        1000 &  5000& 10000 & 20000 \\
  \hline
 $\R_n$ & 0.028 & 0.016 & 0.012 & 0.010 \\
  \hline
  
\end{tabular}
\end{center}

\subsection{Numerical result for non sequential kernel estimator}

Now we give the numerical results for the kernel estimator
defined as

$$
\hat{S}_\zs{n}(z_0)=\dfrac{1}{\sum^n_{k=1}\,Q(u_k)\,y^{2}_\zs{k-1}}\,
\sum^n_{k=1}\,Q(u_k)\,y_{k-1}\,y_\zs{k}\,.
$$
 
For  the standard Gaussian random variables $(\xi_\zs{j})_\zs{j\ge 1}$ in \eqref{sec:In.1},
and by varying the number of observations $n,$ we obtain different risks listed in 
the following table:\\
\begin{center}
\begin{tabular}{|c|c|c|c|c|c|}
  \hline
  $n$ &         1000 &  5000& 10000 & 20000\\
  \hline
 $\R_n$ &  0.046 & 0.026 & 0.02 &0.015  \\
  \hline
  
\end{tabular}
\end{center}

\subsection{Adaptive estimation}

In this case we estimate the function 
\eqref{sec:Nur.1} for 
$z_\zs{0}=1/\sqrt{2}$ and  $\alpha = 0,7$. 
Obviously, that this function belongs to class
$\cH^{(\beta)}(\e,L,L^{*})$ with $\beta=1,7$, $L^{*}=1$ and for any $\e$
and $L$ satisfying the conditions \eqref{sec:Nur.2}.

In the adaptive estimation we take the
 lower regularity 
 $\underline{\beta}=1.6$ and the higher regularity
 $\overline{\beta}=1.8$.

We model the sequential adaptive procedure $\wh{S}_\zs{a,n}=S^{*}_\zs{\check{h}}$  
defined in \eqref{sec:Ad.5}. Numerical results approximate the asymptotic risk 
for this procedure
by the calculation of 
an expectation via $M=30000$ simulations. 
 
$$
\R_\zs{a,n} = \displaystyle\frac{1}{M}\,
\sum_\zs{k=1}^{M}|\wh{S}^{k}_\zs{a,n} - S(z_\zs{0})|
$$                                                                                                                                                                                                                     

By varying the number of observations $n$,
 we obtain different risks listed in the following table:

\begin{center}
\begin{tabular}{|c|c|c|c|c|}
  \hline
  $n$ &        1000 &  5000 &  10000& 20000 \\
  \hline
 $\R_\zs{a,n}$ & 0.021 &0.013  &0.009  &0.007  \\
  \hline
  
\end{tabular}
\end{center}

\medskip

\medskip

\section{Appendix}

\subsection{Concentration properties of the process \eqref{sec:In.1}}

In this section, we study the deviation
\eqref{A.8} for the model \eqref{sec:In.1}.
\begin{lemma}\label{Le.A.1}
For any $q> 1$ and 
$0<\e<1$ the random variables $y_\zs{k}$ in \eqref{sec:In.1} satisfy the following
inequality:
\begin{equation}\label{A.1}
\sup_\zs{n\ge 1}\,\sup_\zs{0\le k\le n}\,
\sup_\zs{S\in\Theta_\zs{\e,L}}\,\sup_\zs{p\in\cP_\zs{\varsigma}}\,
\E_\zs{S,p}\,|y_\zs{k}|^{q}\,\le\,r^{*}_\zs{q}\,, 
\end{equation}
where $r^{*}_\zs{q}$ is defined in \eqref{sec:Pr.1}.
\end{lemma}
\proof From \eqref{sec:In.1} we obtain that for any $k\ge 1$
\begin{equation*}
y_\zs{k}= y_\zs{0}\,\prod_\zs{l=1}^{k}\, S(x_\zs{l})\,+
\sum_\zs{i=1}^{k}\,\prod_\zs{l=i+1}^{k}\, S(x_\zs{l})\,\xi_\zs{i}\,.
\end{equation*} 
Therefore, for $S\in\Theta_\zs{\e,L}$ and  $1\le k\le n$,
$$
|y_\zs{k}|^{q}\le 
2^{q-1}
\left(
|y_\zs{0}|^{q}
+
\left(\sum^k_\zs{j=1}\,(1-\e)^{k-j}\,|\xi_\zs{j}|\right)^{q}
\right)
\,.
$$
Moreover, the H\"older inequality gives
\begin{align*}
\left(\sum^k_\zs{j=1}\,(1-\e)^{k-j}\,|\xi_\zs{j}|\right)^{q}
&\le \left(\sum^k_\zs{j=1}(1-\e)^{k-j}\right)^{q-1}\,\left(\sum^k_\zs{j=1}(1-\e)^{k-j}\,|\xi_\zs{j}|^{q}\right) \\
&\le \left(\frac{1}{\e}\right)^{q-1}\,
\left(\sum^k_\zs{j=1}(1-\e)^{k-j}\,|\xi_\zs{j}|^{q}\right)\,.
\end{align*}
Thus, taking into account the definition \eqref{sec:Sp.0}
we get for any $p\in\cP_\zs{\varsigma}$ 
$$
\E_\zs{S,p}\,
\left(\sum^k_\zs{j=1}\,(1-\e)^{k-j}\,|\xi_\zs{j}|\right)^{q}
\,
\le\,\left(\frac{1}{\e}\right)^{q}\,\s^{*}_\zs{q}\,.
$$
Hence Lemma~\ref{Le.A.1}. \endproof

\noindent
Now we need the following Burkh\"older inequality from \cite{Sh04}.

\begin{lemma}\label{Le.A.2}
Let $(M_\zs{k})_\zs{1\le k\le n}$ be a martingale. Then 
for any $q>1$
\begin{equation}\label{A.1-1}
\E\,|M_\zs{n}|^{q}\le \b^{*}_\zs{q} \E \left(\sum^{n}_\zs{j=1}
(M_\zs{j}-M_\zs{j-1})^{2}\right)^{q/2}\,,
\end{equation}
where the coefficient $\b^{*}_\zs{q}$ is defined in \eqref{sec:Pr.1}.
\end{lemma}

\medskip 

\noindent Now we study the deviation \eqref{A.8}. 

\begin{lemma}\label{Le.A.3}
Let $f$ be a $\bbr\to\bbr$ function twice continuously differentiable 
in $[-1,1]$. Assume also that the bandwidth $h$
satisfies the condition \eqref{sec:Sp.9} -- \eqref{sec:Sp.10}. 
Then for any  $R>0$ and  $q>1$
\begin{equation}\label{A.3} 
\limsup_\zs{n\to \infty}\,
\sup_\zs{k_\zs{*}\le k < m \le k^{*} }\,
\sup_\zs{\underline{\beta}\le \beta \le \overline{\beta}}\,
\sup_\zs{R>0}\sup_\zs{\|f\|_\zs{1}\le R}\,\frac{1}{(Rh)^{q}}
\, \sup_\zs{S \in \Theta_\zs{\e,L}}\,\sup_\zs{p\in\cP_\zs{\varsigma}} \E_\zs{S,p}\,
\left|\varrho_\zs{k,m}(f)\right|^{q}\,\leq \varrho^{*}_\zs{q}\,,
\end{equation}
where
$\|f\|_\zs{1} = \|f\| + \|\dot{f}\|$ and 
$\varrho^{*}_\zs{q}$ is defined in \eqref{sec:Pr.1}.
\end{lemma}
\proof
First of all, note that
\begin{equation}\label{A.4} 
\sum^{m}_\zs{j=k+1}\,f(u_\zs{j})y^2_\zs{j-1}=
T_\zs{k,m}+a_\zs{k,m}\,,
\end{equation}
where 
$T_\zs{k,m}=\sum^{m}_\zs{j=k+1}\,f(u_\zs{j}) y^2_\zs{j}$  and 
$$
a_\zs{k,m}=\sum^{m}_\zs{j=k+1}\,(f(u_\zs{j})-f(u_\zs{j-1}))\,y^2_\zs{j-1}
+f(u_\zs{k})\,y^{2}_\zs{k}
- f(u_\zs{m})\,y^{2}_\zs{m}\,.
$$
Moreover, from 
the model \eqref{sec:In.1} we find 
$$
(1-S^2(z_\zs{0}))T_\zs{k,m}= M_\zs{k,m} + \check{a}_\zs{k,m}+  \sum_\zs{j=k+1}^{m}\,f(u_\zs{j}) 
$$
where $
M_\zs{k,m} =\sum^{m}_\zs{j=k+1}\,(2S(x_\zs{j})\,y_\zs{j-1}\,\xi_\zs{j}\,+
\xi^{2}_\zs{j}-1)\,f(u_\zs{j})$ and
\begin{equation}\label{sec:A.4}
\check{a}_\zs{k,m}=
\sum^{m}_\zs{j=k+1}\,f(u_\zs{j}) S^2(x_\zs{j})y^2_\zs{j-1} - S^2(z_\zs{0})T_\zs{k,m}\,.
\end{equation}
Then we can write $ \varrho_\zs{k,m}(f) $ as follow  
$$
\varrho_\zs{k,m}(f) = \frac{1}{nh\,\gamma(S)}\left(M_\zs{k,m}   + \check{a}_\zs{k,m}\right) + \frac{a_\zs{k,m}}{nh}.
$$
and
\begin{equation}\label{sec:A.6}
\E_\zs{S,p}\,
|\varrho_\zs{k,m}(f)|^{q}\le \frac{3^{q-1}}{\varepsilon^{2q}}\, \E_\zs{S,p}\,
\left(
\left(\frac{|M_\zs{k,m}|}{nh}\right)^{q}
+
\left(\frac{|\check{a}_\zs{k,m}|}{nh}\right)^{q}
+ 
\left(\frac{|a_\zs{k,m}|}{nh}\right)^{q}
\right)\,.
\end{equation}
where $\wt{\kappa}=(1+\kappa_\zs{*})/\kappa_\zs{*}$.
Now we note, that in view of the first condition in \eqref{sec:Sp.10}
for sufficient large $n$
$$
n\kappa_\zs{n}\ge \wt{\kappa}\,,
$$
where $\wt{\kappa}=(1+\kappa_\zs{*})/\kappa_\zs{*}$.
Therefore, for sufficiently large $n$ we get
\begin{equation}\label{sec:A.6-1}
\frac{1}{nh}\le \left(\wt{\kappa}\right)^{2/(2\beta+1)}\,h
\le \wt{\kappa} h
\,,
\end{equation}
\noindent Furthermore, note that 
$(M_\zs{k,j})_\zs{k<j\le m}$ is a martingale.
So, by applying the Burkh\"older inequality \eqref{A.1-1}
and, taking into account that $k^{*}-k_\zs{*} \le 2n h$,
we get
\begin{align*}
\E_\zs{S,p}\,\left(\frac{1}{nh}\,M_\zs{k,m}\right)^{q}
& \le \frac{\b^{*}_\zs{q}R^{q}}{(n h)^{q}}\,\E_\zs{S,p}\,
\left(\sum^{k^{*}}_\zs{j=k_\zs{*}+1}\,
\left(
2\,S(x_\zs{j})\,y_\zs{j-1}\,\xi_\zs{j}+\xi^{2}_j-1\right)^2
\right)^{q/2}\\[2mm]
& 
\le \frac{\b^{*}_\zs{q} R^q}{(nh)^{q/2+1}}\,
\sum^{k^{*}}_\zs{j=k_\zs{*}+1}\E_\zs{S,p}\, \left(
2\,S(x_\zs{j})\,y_\zs{j-1}\,\xi_\zs{j}\,+\,\xi^{2}_\zs{j}-1\right)^{q}\\[2mm]
&\le \frac{4^{q}\b^{*}_\zs{q} R^{q}}{(nh)^{q/2}}\,
\left(r^{*}_\zs{q}\,\s^{*}_\zs{q}+\s^{*}_\zs{2q}+1\right)\le
\,4^{q}\b^{*}_\zs{q}\wt{\kappa}^{q}\, R^{q}\,
\left(r^{*}_\zs{q}\,\s^{*}_\zs{q}+\s^{*}_\zs{2q}+1\right)\,h^{q}
\,,
\end{align*}
\noindent where the coefficients 
$r^{*}_\zs{q}$ and  $\s^{*}_\zs{q}$ are given in \eqref{A.1} and 
\eqref{sec:Sp.0}. 
Note that the term \eqref{sec:A.4} can be rewritten as
\begin{align*}
\check{a}_\zs{k,m}
&=
S^2(z_\zs{0})\left(\sum^{m}_\zs{j=k+1}\,(f(u_\zs{j})- f(u_\zs{j-1}))
\,y^2_\zs{j-1}
+ f(u_\zs{k}) y^{2}_\zs{k}
-f(u_\zs{m}) y^{2}_\zs{m}
\right)\\[2mm]
&+\sum^{m}_\zs{j=k+1}\,f(u_\zs{j}) 
(S^2(x_\zs{j})- S^2(z_\zs{0}))y^2_\zs{j-1}\,.
\end{align*}

\noindent
We recall, that the function $f$ and its derivative $\dot{f}$ are bounded by $R$.
Therefore, taking into account that for all $S\in \Theta_\zs{\e,L}$
 and $k_\zs{*}\le j\le  k^{*}$
the deviation
 $|S(x_\zs{j})-S(z_\zs{0})| \le L|x_\zs{j}-z_\zs{0}|\le Lh$, we obtain
\begin{align*} 
|\check{a}_\zs{k,m}| 
&\leq  R\,
\left(
\left(\frac{1}{nh}+Lh\right)\, 
\sum_\zs{j=k+1}^{m}\,y^{2}_\zs{j}
+y^{2}_\zs{k}+y^{2}_\zs{m}
\right)\\[2mm]
&\le \wt{\kappa}\, R\,
\left((L+1)\,h\,\sum_\zs{j=k+1}^{m}\,y^{2}_\zs{j}
+y^{2}_\zs{k}+y^{2}_\zs{m}
\right)
\,.
\end{align*} 

\noindent 
Therefore, 
\begin{align*}
\sup_\zs{S\in \Theta_\zs{\e,L}}\,
\sup_\zs{p\in\cP_\zs{\varsigma}}\,
\E_\zs{S,p} \left(\frac{1}{nh}\,|\check{a}_\zs{k,m}|\right)^{q}
&\le 2^{q-1}\,\wt{\kappa}^{q}
R^{q}r^{*}_\zs{2q}
\left(
2^{q}(1+L)^{q}h^{q}
+\frac{2^{q}}{(nh)^{q}}
\right)\\[2mm]
&\le 
\left(4R(1+L) 
\right)^{q}r^{*}_\zs{2q} h^{q}\,.
\end{align*}
\noindent 
Similarly,
$$
\sup_\zs{S\in \Theta_\zs{\e,L}}\,
\sup_\zs{p\in\cP_\zs{\varsigma}}\,
\E_\zs{S,p} \left(\frac{1}{nh}\,|a_\zs{k,m}|\right)^{q}
\le 3\,2^{q-1}\,\wt{\kappa}^{q}\,R^{q}\,r^{*}_\zs{2q}\, h^{q}\,.
$$
Then, taking this into account in \eqref{sec:A.6} 
we obtain the upper bound \eqref{A.3}. Hence Lemma~\ref{Le.A.3}. 
\endproof

\subsection{Uniform limit theorem}

In this section we study the following sequence 
\begin{equation}\label{sec:A.7}
\wt{\zeta}_\zs{n}(h)=
\zeta_\zs{n}(h)\,\Chi_\zs{(A_\zs{\nu,n}\geq H)}\,,
\end{equation}
where $\zeta_\zs{n}(h)$ defined by \eqref{sec:Prf.3}, the bandwidth $h$ 
is defined \eqref{sec:Sp.9} and the threshold $H$ is given in  \eqref{sec:Sp.8}.
We will make use of the following result.
\begin{lemma}\label{Le.A.4}(cf. \cite{Fr71}, p. 90-91)
 Let $0<\delta<1$ and $r>0$. Assume that
 $\left(\m_\zs{k}\right)_\zs{k\ge 1}$ is a martingale difference with respect to 
 the filtration $(\cF_\zs{k})_\zs{k\ge 1}$ such that
$$
|\m_\zs{k}|\le \delta r^{1/2}
\quad\mbox{and}\quad
\sum^{\infty}_\zs{k=1}\,\E\left(\m^{2}_\zs{k}|\cF_\zs{k-1}\right)\ge r\,.
$$
Let
$$
\tau=\inf\left\{k\ge 1: \sum^{k}_\zs{j=1}\,\E(\m^{2}_\zs{j}|\cF_\zs{j-1})\ge r
\right\}\,.
$$
There exists a function $\rho\,:\,(0,+\infty)\to [0,2]$ not depending
on distribution of the martingale difference, such that
$\lim_\zs{x\to 0}\,\rho(x)=0$ and
$$
\sup_\zs{x\in \R}\,
\left|
\P\left(
\frac{1}{r^{1/2}}\sum^{\tau}_\zs{k=1}\,\m_\zs{k}\,\le x\right) 
- 
\Phi(x)
\right|\le 
\rho(\delta)\,, 
$$ 
where $\Phi$ is the standard normal distribution function.
\end{lemma}

\begin{lemma}\label{Le.A.5}
The sequence \eqref{sec:Prf.3} satisfies the following
limiting property:
 $$
\wt{\zeta}_\zs{n}\Longrightarrow\,\zeta \sim {\cal N}(0,1)
\quad\mbox{uniformly in}
\  p\in \cP_\zs{\varsigma}
\quad\mbox{and}\quad
S\in\Theta_\zs{\e,L}\,.
$$
\end{lemma}
\proof First  for some $0<\delta<1$ we set
$$
\m_\zs{j}=
\,Q(u_\zs{j})\,\check{y}_\zs{j-1}\,
\check{\xi}_\zs{j}\Chi_\zs{\{\nu<j\le n\}}+
\delta \check{\xi}_\zs{j}\Chi_\zs{\{j>n\}}\,,
$$
where 
$\check{y}_\zs{j}=y_\zs{j}\,\Chi_\zs{\{|y_\zs{j}|\le \delta^{2}\,\check{H}^{1/2}\}}$,
$\check{H}=\E_\zs{p}\,\check{\xi}^{2}_\zs{1}\, H$
and 
 $\check{\xi}_\zs{j}=\xi_\zs{j}\,\Chi_\zs{\{|\xi_\zs{j}|\le \delta^{-1}\}}-
\E_\zs{p}\,\xi_1\,\Chi_\zs{\{|\xi_1|\le \delta^{-1}\}}$. It is clear that the
sequence
$(\m_\zs{\nu+j})_\zs{j\ge 0}$ is a martingale difference  with respect to
$(\cG_\zs{j})_\zs{j\ge 0}$, where $\cG_\zs{j}$ is $\sigma$ - field generated by 
the observations $\{y_\zs{1}\,,\ldots,y_\zs{\nu+j}\}$.
Now we set
$$
\check{\zeta}_\zs{\check{H}}=\frac{1}{\sqrt{\check{H}}}\,\sum^{\check{\tau}_\zs{H}}_\zs{j=1}\,
\m_\zs{\nu+j}\,,
$$
where
$\check{\tau}_\zs{H}=\inf\left\{k\ge 1\,:\,
 \sum^{k}_\zs{j=1}\,\E\,(\m^{2}_\zs{\nu+j}|\cG_\zs{j-1})\ge \check{H}\right\}$.
Note, that $\check{\tau}_\zs{H}=\tau_\zs{H}$ on the set $\{A_\zs{\nu,n}\geq H\}$
 for any $H>0$. Lemma~\ref{Le.A.4} implies that
$\check{\zeta}_\zs{\check{H}}$  goes in distribution to $\cN(0,1)$ uniformly
in $p\in \cP_\zs{\varsigma}$ and $S\in\Theta_\zs{\e,L}$ as 
$\delta\to 0$. Now we set
$$
\check{\Omega}_\zs{n}=\cap^{k^{*}}_\zs{j=\nu}\,\{y_\zs{j}= \wt{y}_\zs{j}\}\,.
$$
Using Lemma~\ref{Le.A.1} through the Chebyshev inequality we obtain that
$$
\sup_\zs{S\in\Theta_\zs{\e,L}}\,
\sup_\zs{p\in \cP_\zs{\varsigma}}\,
\P_\zs{S,p}\left( \check{\Omega}^{c}_\zs{n}\right)\le \,
\frac{(k^{*}-\nu)\,r^{*}_\zs{4}}{\delta^{8} \E_\zs{p}\check{\xi}^{2}_\zs{1} (nh)^{2}}
\to 0
\quad\mbox{as}\quad n\to\infty\,.
$$
Moreover, note that on the set $\check{\Omega}_\zs{n}\cap \{A_\zs{\nu,n}\geq H\}$
\begin{equation}\label{sec:A.8}
\left(\frac{\check{H}}{H}\right)^{1/2}\wt{\zeta}_\zs{n}
-\check{\zeta}_\zs{\check{H}}
=
\check{\Delta}_\zs{1}
+
\check{\Delta}_\zs{2}
\,,
\end{equation}
where
$$
\check{\Delta}_\zs{1}=
\frac{1}{\sqrt{\check{H}}}\sum^{\tau_\zs{H}}_\zs{j=\nu+1}(\varkappa_\zs{j}-1)\,
Q(u_\zs{j})\,\check{y}_\zs{j-1}\,\xi_\zs{j}\,,
\quad
\check{\Delta}_\zs{2}=
\frac{1}{\sqrt{\check{H}}}\sum^{\tau_\zs{H}}_\zs{j=\nu+1}\,
Q(u_\zs{j})\,\check{y}_\zs{j-1}\,\wt{\xi}_\zs{j}
$$
and  $\wt{\xi}_\zs{j}=\xi_\zs{j}-\check{\xi}_\zs{j}=
\xi_\zs{j}\,\Chi_\zs{\{|\xi_\zs{j}|> \delta^{-1}\}}-
\E\,\xi_1\,\Chi_\zs{\{|\xi_1|> \delta^{-1}\}}$.
Note, that
$$
\E_\zs{S,p}\left( 
\check{\Delta}^{2}_\zs{1}
|\cG_\zs{0}
\right)
\le \delta^{4} \E_\zs{S,p}\left(
\sum^{k^{*}}_\zs{j=\nu+1}(1-\varkappa_\zs{j})^{2}
|\cG_\zs{0}
\right)\le \delta^{4}\,.
$$
Moreover, taking into account that $\check{y}^{2}_\zs{j}\le y^{2}_\zs{j}$,
we get
\begin{align*}
\E_\zs{S,p}\left( 
\check{\Delta}^{2}_\zs{2}
|\cG_\zs{0}
\right)
&\le
\frac{\E_\zs{p}\,\wt{\xi}^{2}_\zs{1}}{\check{H}}
 \E_\zs{S,p}\left(
\sum^{\tau_\zs{H}}_\zs{j=\nu+1}\,
Q(u_\zs{j})\,\check{y}^{2}_\zs{j-1}
|\cG_\zs{0}
\right)\\[2mm]
&\le
\frac{\E_\zs{p}\,\wt{\xi}^{2}_\zs{1}}{\check{H}}
\left(
H+\delta^{4} \check{H}
\right)=
\E_\zs{p}\,\wt{\xi}^{2}_\zs{1}
\left(
\frac{1}{\E_\zs{p}\,\check{\xi}^{2}_\zs{1}}
+\delta^{4}
\right)
\end{align*}
Taking into account here, that 
$$
\lim_\zs{\delta\to 0}\sup_\zs{p\in\cP_\zs{\varsigma}}\,\E_\zs{p}\,\wt{\xi}^{2}_\zs{1}\,=0
\quad\mbox{and}\quad
\lim_\zs{\delta\to 0}\sup_\zs{p\in\cP_\zs{\varsigma}}
\left|\E_\zs{p}\,
\check{\xi}^{2}_\zs{1}-1\right|=0\,,
$$
we obtain 
$$
\lim_\zs{\delta\to 0}
\sup_\zs{S\in\Theta_\zs{\e,L}}\,
\sup_\zs{p\in \cP_\zs{\varsigma}}\,
\max\left( 
\E_\zs{S,p}\, 
\check{\Delta}^{2}_\zs{1}
\,,\,
\E_\zs{S,p}\, 
\check{\Delta}^{2}_\zs{2}
\right)
=0\,.
$$
Therefore,  Proposition~\ref{Pr.S.3} 
and the representation  \eqref{sec:A.8} yield  for any $\mu>0$
$$
\lim_\zs{\delta\to 0}\,
\limsup_\zs{n\to\infty}\,
\sup_\zs{S\in\Theta_\zs{\e,L}}\,
\sup_\zs{p\in \cP_\zs{\varsigma}}\,
\P_\zs{S,p}\left(
\left|\wt{\zeta}_\zs{n}
-\check{\zeta}_\zs{\check{H}}\right|>\mu
\right)=0\,.
$$
Hence Lemma~\ref{Le.A.5}. \endproof

\medskip

\subsection{Properties of $\wt{\zeta}_\zs{n}(h)$}

\begin{lemma}\label{Le.A.6} 
For all $z\ge 2$ 
\begin{equation}\label{sec:A.9}
\sup_\zs{n\ge 1}
\sup_\zs{h>0}\,\sup_\zs{S\in \C[0,1]}\,
\P_\zs{S}\left( 
\wt{\zeta}_\zs{n}(h)\ge z
\right)
\le 2 e^{-z^{2}/8}\,.
\end{equation}
\end{lemma}
\noindent 
The proof of this Lemma is the same as Lemma A.5 from \cite{Ark11}.

\noindent 
Using this lemma we can obtain that for any $q>2$

\begin{equation}\label{sec:A.10}
\sup_\zs{n\ge 1}
\sup_\zs{h>0}\,\sup_\zs{S\in \C[0,1]}\,
\E_\zs{S}\,
|\wt{\zeta}_\zs{n}(h)|^{q}
\,
\le \mu^{*}_\zs{q}\,,
\end{equation}
where $\mu^{q}_\zs{q}=2^{q}+2q\int^{\infty}_\zs{0}t^{q-1}\,e^{-t^{2}/2}\d t$.

\medskip

\medskip
\medskip

\bibliographystyle{plain}


\end{document}